\documentclass{article}
\usepackage{amsmath,amsxtra,amssymb,latexsym,amsfonts,amscd}
\usepackage{multicol,color}
\usepackage{float}
\usepackage{soul}
\usepackage{graphicx}
\usepackage{pstricks-add}
\usepackage{pgf,tikz}
\usepackage{mathrsfs}
\usetikzlibrary{arrows}
\usepackage{amssymb}
\usepackage{theorem}
\usepackage{fancyhdr}
\usepackage{tikz}
\usepackage{enumerate}
\usepackage[margin=1.1in]{geometry}
\usepackage[normalem]{ulem}
\usepackage{bbm}
\def\disp{\displaystyle}
\def\e{\varepsilon}
\def\ve{\varepsilon}
\def\dd{\delta}
\def\DD{\Delta}
\def\lm{\lambda}
\def\O{\Omega}

\def\({\left(}
\def\){\right)}
\def\[{\left[}
\def\]{\right]}
\def\n{\left \|}
\def\en{\right \|}
\def\nn{\left \{ }
\def\hnn{\right \}}

\def\ox{\bar{x}}
\def\oy{\bar{y}}
\def\oz{\bar{z}}
\def\ov{\bar{v}}
\def\ou{\bar{u}}

\def\gph{\hbox{}}

\def\gg{\gamma}
\def\dn{\downarrow}

\def\tto{\rightrightarrows}

\def\Bar{\overline}
\def\la{\langle}
\def\ra{\rangle}
\def\ve{\varepsilon}
\def\B{I\!\!B}
\def\h{\hfill\Box}
\def\R{\mathbb{R}}
\def\N{\mathbb{N}}

\def\co{\mbox{\rm co}\,}
\def\gph{\mbox{\rm gph}\,}
\def\epi{\mbox{\rm epi}\,}

\def\dom{\mbox{\rm dom}\,}

\def\dist{\mbox{\rm dist}}

\def\dn{\downarrow}
\def\O{\Omega}
\def\o{\omega}

\def\vph{\varphi}
\def\emp{\emptyset}

\def\oR{\Bar{\R}}
\def\lm{\lambda}

\def\gg{\gamma}

\def\dd{\delta}
\def\DD{\Delta}
\def\al{\alpha}
\def\vth{\vartheta}

\def\vth{\vartheta}
\def\ph{\varphi}

\def\N{I\!\!N}
\def\th{\theta}

\def\ne{\neq}
\newtheorem{theorem}{Theorem}[section]

\newtheorem{definition}[theorem]{Definition}
\theoremstyle{plain}{\theorembodyfont{\rmfamily}
}
\theoremstyle{plain}{\theorembodyfont{\rmfamily}
}
\theoremstyle{plain}{\theorembodyfont{\rmfamily}
}
\theoremstyle{plain}{\theorembodyfont{\rmfamily}
}

\theoremstyle{plain}{\theorembodyfont{\rmfamily}
}

\def\eq{\begin{equation}}
\def\eeq{\end{equation}}
\begin{document}
\begin{center}
{\bf DISCRETE APPROXIMATIONS AND OPTIMAL CONTROL\\ OF NONSMOOTH PERTURBED SWEEPING PROCESSES}\footnote{This research was partly supported by the USA National Science Foundation under grants DMS-1512846 and DMS-1808978 and the USA Air Force Office of Scientific Research grant \#15RT0462}\\[3ex]
BORIS S. MORDUKHOVICH\footnote{Department of Mathematics, Wayne State University, Detroit, Michigan 48202, USA (boris@math.wayne.edu). Research of this author was also supported in part by the Australian Research Council under Discovery Project DP-190100555.}\;\;and\;\;DAO NGUYEN\footnote{Department of Mathematics, Wayne State University, Detroit, Michigan 48202, USA (dao.nguyen2@wayne.edu).}\\[2ex]
{\bf Dedicated to Umberto Mosco in honor of his 80th birthday}\\[1ex]
\end{center}
\small{\sc Abstract.} The main goal of this paper is developing the method of discrete approximations to derive necessary optimality conditions for a class of constrained sweeping processes with nonsmooth perturbations. Optimal control problems for sweeping processes have been recently recognized among the most interesting and challenging problems in modern control theory for discontinuous differential inclusions with irregular dynamics and implicit state constrained, while deriving necessary optimality conditions for their local minimizers have been significantly based on the smoothness of controlled dynamic perturbations. To overcome these difficulties, we use the method of discrete approximations and employ advanced tools of second-order variational analysis. This approach allows us to obtain new necessary optimality conditions for nonsmooth and nonconvex discrete-time problems of the sweeping type. Then we employ the obtained conditions and the strong convergence of discrete approximations to establish novel results for original nonsmooth sweeping control problems that include extended Euler-Lagrange and maximization conditions for local minimizers. Finally, we present applications of the obtained results to solving a controlled mobile robot model with a nonsmooth sweeping dynamics that is of some practical interest.\\[1ex]
\vspace*{-0.2in}

\normalsize
\section{Introduction and Overview}\label{assumptions}
\setcounter{equation}{0}\vspace*{-0.1in}

The basic sweeping process was introduced by Jean Jacques Moreau in the beginning of 1970s being mainly motivated by applications to elastoplasticity; see \cite{moreau} for more references and discussions. It was described by the dissipative differential inclusion
\begin{equation}\label{SP}
-\dot x(t)\in N\big(x(t);C(t)\big)\;\mbox{ a.e. }\;t\in[0,T]\;\mbox{ with }\;x(0):=x_0\in C(0)
\end{equation}
via the normal cone to a continuously moving convex set $C(t)$, where the normal cone to a convex set is
\begin{equation}\label{NC}
N(x;\O)=N_\O(x):=\big\{v\in\R^n\big|\;\la v,y-x\ra\le 0,\;y\in\O\big\}\;\textrm{if}\;x\in\O\;\textrm{ and }\;N(x;\O):=\emp\textrm{ if }x\notin\O.
\end{equation}
Over the years, mathematical theory of the sweeping process \eqref{SP} and its modifications has been largely developed and applied to many practical models; see, e.g., \cite{bt,CT,HB,KMM,mv,Thi} and the references therein.\vspace*{-0.05in}

One of the most fundamental results of the sweeping process theory for model \eqref{SP} and its generalizations is the existence and {\em uniqueness} of solutions to the Cauchy problem as in \eqref{SP}. This clearly excludes to possibility to optimize the sweeping dynamics. It is mainly due to {\em discontinuity} and {\em monotonicity} properties of the normal cone mapping governed the sweeping dynamics. Note that this issue dramatically differs sweeping processes from {\em Lipschitzian} differential inclusions for which optimal control theory has been comprehensively developed; see the monographs \cite{cl,m-book2,v} among other publications.\vspace*{-0.05in}

{\em Optimal control} problems for sweeping processes were formulated and investigated more recently. Although the existence and relaxation questions for the sweeping dynamics with controlled perturbations appeared earlier \cite{et}, the first necessary optimality conditions for a novel class of problems in dynamic optimization concerning systems \eqref{SP} with controlled moving sets $C(t)=C(u(t))$  were obtained in \cite{chhm1}. Then necessary optimality optimality conditions of different types were derived in the series of publications dealing with several versions of controlled sweeping processes (see, e.g., \cite{ao,ac,bk,cm1,cm2,cm3,chhm3,hm,pfs}). By now optimal control theory for such systems has become a very active and attractive area of research in dynamic optimization with applications to practical modeling.\vspace*{-0.05in}

It has been realized that major difficulties in deriving necessary optimality conditions for controlled sweeping process come not only from the high discontinuity and irregularity of the sweeping dynamics, but also from the intrinsic presence of {\em pointwise state} and {\em mixed state-control constraints} that are challenging and underinvestigated even in the classical control theory for smooth systems; see, e.g., the very recent survey \cite{ak}. To overcome these difficulties, {\em all} the known approaches to deriving necessary optimality conditions for controlled sweeping processes employ one or another {\em approximation procedure}, with an involved convergence analysis. Note that various approximation/perturbation procedures and appropriate convergence notions have always been among the most efficient tools of nonlinear analysis, particularly of its variational aspects. Recall to this end the celebrated notion of {\em Mosco convergence} \cite{mosco} (equivalent to {\em epi-convergence} in finite dimensions), the extended notion of {\em $\Gamma$-convergence} coming from the school of De Giorgi, etc.; see \cite{abm,m18,rw} with the vast bibliographies therein.\vspace*{-0.03in}

Having said this, let us now formulate the problem of our study in this paper and discuss the approximation method of deriving necessary optimality conditions for its local minimizers. The {\em sweeping optimal control} problem $(P)$ is on minimizing the cost functional
\begin{equation}\label{cost1}
J[x,u]:=\vph\big(x(T)\big)
\end{equation}
over {\em feasible solutions} $(x(\cdot),u(\cdot))\in W^{1,2}([0,T];\R^n)\times L^2([0,T];\R^d)$ satisfying the constrained system
\begin{equation}\label{Problem}
\left\{\begin{matrix}
-\dot{x}(t)\in N\big(x(t);C\big)+g\big(x(t),u(t)\big)\;\textrm{ a.e. }\;t\in[0,T],\;x(0)=x_0\in C\subset\R^n,\\
u(t)\in U\subset\R^d\;\textrm{ a.e. }\;t\in[0,T],
\end{matrix}\right.
\end{equation}
where the set $C\subset\R^n$ is a convex polyhedron defined by
\begin{equation}\label{C}
C:=\bigcap_{j=1}^{s}C^j\;\textrm{ with }\;C^j:=\nn x\in\R^n\big|\;\la x^j_*,x\ra\le c_j\hnn.
\end{equation}
It follows from the normal cone definition \eqref{NC} and the sweeping differential inclusion in \eqref{Problem} that any feasible trajectory $x(\cdot)$ automatically satisfies the {\em pointwise state constraints}
\begin{equation}\label{state}
\la x^j_*,x(t)\ra\le c_j\;\mbox{ for all }\;t\in[0,T]\;\mbox{ and }\;j=1,\ldots,s.
\end{equation}\vspace*{-0.2in}

Necessary optimality conditions for optimal control problems governed by sweeping processes with controlled perturbations were first obtained in \cite{cm1} with control actions entering not only additive perturbations, but also the moving sweeping set $C(\cdot)$. However, the problems considered in \cite{cm1} and in the subsequent papers \cite{cm2,cm3} (the latter dealt with the case of nonconvex moving sets) did not impose any constraints on feasible controls $u(\cdot)$, which were assumed to be smooth, namely $u(\cdot)\in W^{1,2}([0,T];\R^d)$.\vspace*{-0.05in}

More conventional---from the viewpoint of basic control theory---sweeping controls models of type \eqref{Problem} with pointwise constraints on measurable controls were studied quite recently in \cite{ac,cmn1,cmn,pfs} under different assumptions on the sweeping set $C$. A common feature of these and other papers on perturbed controlled sweeping processes was the {\em smoothness} assumption on the perturbation function $g(x,u)$ with respect to $x$ in \cite{pfs} and with respect to both variables $(x,u)$ in \cite{ac,cm1,cm2,cm3,cmn1,cmn}. Furthermore, in \cite{cmn1,cmn} the partial Jacobian $\nabla_u g$ was assumed to have full rank at the optimal solution. The smoothness and full rank assumptions imposed in those papers were dictated by and essentially used in the approximation techniques employed in the derivation of the corresponding optimality conditions.\vspace*{-0.03in}

The sweeping optimal control model investigated in our paper agrees with those in \cite{cmn1,cmn}, while we avoid here the aforementioned {\em smoothness} and {\em full rank assumptions} to derive necessary optimality conditions for local minimizers of $(P)$. Similarly to \cite{cmn1,cmn}, our approach is based on the method of {\em discrete approximations}, but the main novelty consists of using advanced robust tools and results of {\em second-order variational analysis} and {\em generalized differentiation} to deal with nonsmooth mappings and to justify an adequate convergence procedure. In this way we obtain new necessary optimality conditions for nonsmooth and nonconvex control problems of discrete approximation that lead us by passing to the limit to novel relationships for relaxed local minimizers of the original sweeping control problems with totally nonsmooth data. The given application to a controlled mobile robot model with nonsmooth dynamics confirms the efficiency of the obtained results.\vspace*{-0.03in}

The rest of the paper is organized as follows. In Section~\ref{assump} we formulate and discuss the concept of intermediate local minimizers of $(P)$ for which necessary optimality conditions is derived below. Section~\ref{tools} recalls the basic constructions of first-order generalized differentiation for sets, set-valued mappings, and extended-real-valued functions that are used in the paper. Section~\ref{2nd-order} concerns second-order constructions of variational analysis and presents their calculations in terms of the given data of the sweeping control problem $(P)$, which are crucial for the subsequent material. In Section~\ref{disc} we construct a well-posed sequence of discrete approximations and establish necessary optimality conditions for their optimal solutions, which strongly converge to a designated local minimizer of $(P)$. Section~\ref{optim} contains the main results on necessary optimality conditions for the local minimizes of $(P)$ under investigation obtained by passing to the limit from those derived in Section~\ref{disc}. The concluding Section~\ref{applic} presents applications of the obtained necessary optimality conditions to solving a controlled version of a mobile robot model of some practical interest that is described by a nonsmooth perturbed sweeping process.\vspace*{-0.2in}

\section{Local Minimizers and Standing Assumptions}\label{assump}
\setcounter{equation}{0}\vspace*{-0.1in}

We begin with defining the notion of local minimizers for the sweeping control problem $(P)$.\vspace*{-0.1in}

\begin{definition}[\bf local minimizers]\label{ilm} Let $(\ox(\cdot),\ou(\cdot))$ be a feasible pair for the control problem $(P)$. It is said to be a $W^{1,2}\times L^2$-{\sc local minimizer} for $(P)$ if there exists $\e>0$ such that $J[\ox,\ou]\le J[x,u]$ for any pair $(x(\cdot),u(\cdot))$, which is feasible to $(P)$ and satisfies the localization condition
\begin{equation}\label{local}
\int_0^T\(\n\dot x(t)-\dot{\ox}(t)\en^2+\n u(t)-\ou(t)\en^2\)dt<\e.
\end{equation}
\end{definition}\vspace*{-0.1in}

For differential inclusions in the form $\dot x\in F(t,x)$, with no explicit control term, local minimizers from Definition~\ref{ilm} reduce to {\em intermediate local minimizers} of rank two introduced in \cite{m95}. The name comes from the fact this notion occupied an intermediate position between the classical concepts of weak and strong minimizers in variational problems. In the books \cite{m-book2,v} and the references therein the reader can find more information, examples, and results on intermediate local minimizers for Lipschitzian differential inclusions that play an independent role for such problems even in simple settings; cf.\ also \cite{cm3,chhm3,cmn,hm} for controlled sweeping processes. If the space $W^{1,2}[0,T]$ is replaced by ${\cal C}[0,T]$ in Definition~\ref{ilm}, i.e., the the norm $\int_0^T\|\dot{x}(t)-\dot{\ox}(t)\|^2dt$ is replaced by $\max_{t\in[0,T]}\|x(t)-\ox(t)\|$ in \eqref{local}, then we arrive at the notion of {\em strong} ${\cal C}\times L^2$-local minimizers for $(P)$. The latter minimizers are clearly included in the collection of $W^{1,2}\times L^2$ ones from Definition~\ref{ilm}, and the inclusion is generally strict.\vspace*{-0.03in}

Let us now formulate the {\em standing assumptions} imposed throughout the paper without further mentioning. Some of them involve a reference trajectory $\ox(\cdot)$, which is specified below as (a part of) either a given local minimizer for $(P)$, or an optimal solution to its discrete approximation.\\[1ex]
{\bf(H1)} The cost function $\ph\colon\R^n\to\R$ in \eqref{cost1} is {\em locally Lipschitzian} around $\ox(T)$.\\
{\bf(H2)} The control region $U\ne\emp$ is a {\em compact} set in $\R^d$.\\
{\bf(H3)} The perturbation mapping $g\colon\R^n\times\R^d\to\R^n$ from \eqref{SP} is {\em Lipschitz continuous} with respect to both variables $x$ and $u$ uniformly on $U$ whenever $x$ belongs to a bounded subset of $\R^n$, and there exists a number $\beta>0$ for which the following {\em sublinear growth condition} holds:
\begin{equation*}
\|g(x,u)\|\le\beta\big(1+\|x\|\big)\;\mbox{ for all }\;u\in U.
\end{equation*}
{\bf(H4)} For the convex polyhedron $C$ from \eqref{C} with the vertices $x^j_*$, $j=1,\ldots,s$, the {\em positive linear independence constraint qualification} (PLICQ) condition
\begin{equation}\label{plicq}
\Big[\sum_{j\in I(\ox)}\al_j x^j_*=0,\;\al_j\ge 0\Big]\Longrightarrow\big[\al_j=0\;\textrm{ for all }\;j\in I(\ox)\big\}
\end{equation}
is satisfied along $\ox=\ox(t)$ for all $t\in[0,T]$, where the set of {\em active constraint indices} is given by
\begin{equation}\label{aci}
I(\ox):=\big\{j\in\{1,\ldots,s\}\;\big|\;\la x^j_*,\ox\ra=c_j\big\}.
\end{equation}\vspace*{-0.2in}

Note that for polyhedral under consideration the PLICQ condition reduces to the {\em Mangasarian-Fromovitz constraint qualification}, which is a major qualification condition in nonlinear programming. This condition is much less restrictive than the {\em linear independent constraint qualification} (LICQ), which corresponds to \eqref{plicq} with the replacement of $\al_j\ge 0$ on the left-hand side by $\al_j\in\R$ for $j\in I(\ox)$.\vspace*{-0.03in}

To proceed further, we need to slightly modify the notion of local minimizers for which the necessary optimality conditions are derived in what follows. It has been well recognized (starting with the pioneering work by Bogolyubov and Young on the classical calculus of variations in the 1930s) that the study of variational and control problems involving time derivatives requires a certain {\em convexification} with respect to derivative/velocity variables in order to conduct limiting procedures. Furthermore, such a convexified extension (or {\em relaxation} in the terminology coined by Warga in 1962) ensures the existence of relaxed optimal solutions and often keeps the {\em same value} of the cost functional. For conventional optimal control problems and Lipschitzian differential inclusions, the reader can find more results, references, and discussions in \cite{m-book2,v}. We also refer the reader to \cite{dfm,et,Tol} to relaxation procedures and results of the aforementioned types for non-Lipschitzian differential inclusions and controlled sweeping processes.\vspace*{-0.03in}

Following this line, we define now the needed relaxation of the $W^{1,2}\times L^2$-local minimizers for problem $(P)$ under the standing assumptions. Performing the convexification of the differential inclusion in \eqref{SP} with taking into account the convexity of $N(x;C)$ and compactness of the image sets $g(x,U)$ by (H2) and (H3), we arrive at the convexified inclusion
\begin{equation}\label{co}
-\dot{x}(t)\in N\big(x(t);C\big)+\co g\big(x(t),U\big)\;\textrm{ a.e. }\;t\in[0,T],\;x(0)=x_0\in C,
\end{equation}
where ``co" stands for the convex hull of a set.\vspace*{-0.1in}

\begin{definition}{\bf(relaxed local minimizers).}\label{relaxed} Given a pair $(\ox(\cdot),\ou(\cdot))$ feasible to $(P)$, we say that it is a {\sc relaxed $W^{1,2}\times L^2$-local minimizer} for the original problem if there exists $\e>0$ such that
\begin{equation*}
\ph\big(\ox(T)\big)\le\ph\big(x(T)\big)\;\textrm{ whenever }\;\int_0^T\(\n\dot x(t)-\dot{\ox}(t)\en^2+\n u(t)-\ou(t)\en^2\)dt<\e,
\end{equation*}
where $u(\cdot)$ is a measurable control with $u(t)\in\co U$ a.e.\ on $[0,T]$, and where $x(\cdot)$ is a trajectory of the convexified inclusion \eqref{co} that can be strongly approximated in $W^{1,2}([0,T];\R^n)$ by feasible trajectories to $(P)$ generated by piecewise constant controls $u_m(\cdot)$ on $[0,T]$ the convex combinations of which strongly converge to $u(\cdot)$ in the norm topology of $L^2([0,T];\R^d)$.
\end{definition}\vspace*{-0.05in}

It is easy to see that there is no difference between $W^{1,2}\times L^2$-local minimizers of $(P)$ and their relaxation in Definition~\ref{relaxed} if the sets $g(x,U)$ are {\em convex}, which is not assumed in what follows. Furthermore, a close look at the proof of \cite[Theorem~2]{et} allows us to deduce that any {\em strong} ${\cal C}\times L^2$-local minimizers for $(P)$ is a relaxed one under the imposed standing assumptions.\vspace*{-0.2in}

\section{Preliminary from Generalized Differentiation}\label{tools}
\setcounter{equation}{0}\vspace*{-0.1in}

In this section we briefly review some basic constructions of first-order generalized differentiation in variational analysis that are employed in the paper. Note that, although the set under consideration in the original sweeping model \eqref{SP} is convex, we have to deal with generalized normals to nonconvex sets and with the corresponding first-order and second-order constructions for extended-real-valued functions and set-valued mappings. This is due to the essence of the method of discrete approximations, which reduces optimization of differential inclusions to mathematical programs having  many geometric constraints of the (nonconvex) {\em graphical} type with the subsequent passing to the limit. Moreover, the generalized differential constructions that are suitable for such a device should be robust, have rich calculus, and be able to adequately handle graphical sets. These requirements are satisfied for generalized differentiation theory initiated by the first author; see \cite{m-book1,m18,rw} for more details and references.\vspace*{-0.03in}

Given a set-valued mapping $F\colon\R^n\tto\R^m$ and a point $\ox$ with $F(\ox)\ne 0$, the symbol
\begin{equation*}
\underset{x\to\ox}{\textrm{Lim sup }}F(x):=\{z\in\R^m\left|\right.\exists\textrm{ sequences }x_k\to\ox,\,z_k\to z\textrm{ such that }z_k\in F(x_k),\,k=1,2,\ldots\}
\end{equation*}
denotes the (Kuratowski-Painlev\'e) $\textit{outer limit}$ of $F$ at $\ox$. For a set $\O\subset\R^n$ locally closed around $\ox\in\O$, the (Mordukhovich basic/limiting) {\em normal cone} to $\O$ at $\ox$ defined by
\begin{equation}\label{nor_con}
N(\ox;\O)=N_{\O}(\ox):=\underset{x\to\ox}{\textrm{Lim sup}}\big\{\textrm{cone}\big[x-\Pi(x;\O)\big]\big\},
\end{equation}
where $\Pi(\ox;\O)$ stands for the Euclidean projection of $\ox$ onto $\O$ given by
$$
\Pi(\ox;\O):=\big\{y\in\O\,\big{|}\,\|\ox-y\|=\dist(\ox,\O)\big\}.
$$
When $\O$ is convex, the normal cone \eqref{nor_con} reduces to a classical one of convex analysis \eqref{NC}, while in general the cone is nonconvex. However, this normal cone together with the associated subdifferential and coderivative constructions enjoys {\em full calculus} based on {\em variational/extremal principles}; see \cite{m-book1,rw}. Note that the {\em convex closure} of \eqref{nor_con} agrees in finite dimensions with the normal cone introduced by Clarke via the duality correspondence with his tangent cone to $\O$ at $\ox$: see \cite{cl,m-book1,rw} for more details and references. For one of the results of Theorem~\ref{nsThm6.1*} below, we use the convex hull $\co N(\ox;\O)$ without taking the closure operation. This allows us to exploit the {\em robustness} (outer semicontinuity) property of the normal cone mapping $x\mapsto N(x;\O)$ from \eqref{nor_con}, which is not shared by the normal cone of Clarke.\vspace*{-0.03in}

Let $F\colon\R^n\tto\R^m$ with its domain and graph given by
\begin{equation*}
\dom F:=\big\{x\in\R^n\;\big|\;F(x)\ne\emp\big\}\;\mbox{ and }\;\gph F:=\big\{(x,y)\in\R^n\times\R^m\;\big|\;y\in F(x)\big\},
\end{equation*}
respectively. Assuming that the graph of $F$ is locally closed around $(\ox,\oy)\in\gph F$, the {\em coderivative} of $F$ at $(\ox,\oy)$ is defined via the normal cone \eqref{nor_con} by
\begin{equation}\label{cod}
D^*F(\ox,\oy)(u):=\big\{v\in\R^n\left|\right.(v,-u)\in N\big((\ox,\oy);\gph F\big)\big\},\quad u\in\R^m.
\end{equation}
When $F\colon\R^n\to\R^m$ is single-valued and ${\cal C}^1$-smooth around $\ox$, we get the representation
$$
D^*F(\ox)(u)=\big\{\nabla F(\ox)^*u\big\}\;\textrm{ for all }\;u\in\R^m
$$
via the adjoint/transposed Jacobian matrix $\nabla F(\ox)^*$, where $\oy=F(\ox)$ is omitted. If $F\colon\R^n\to\R^m$ is single-valued and locally Lipschitzian around $\ox$, then
\begin{equation*}
\co D^*F(\ox)(u)=\big\{A^*u\;\big|\;A\in\Bar\nabla F(\ox)\big\},\quad u\in\R^n,
\end{equation*}
where $\Bar\nabla F(\ox)$ is the (Clarke) {\em generalized Jacobian} of $F$ at $\ox$, which a nonempty compact subset of the matrix spaces $\R^{m\times n}$ defined as
\begin{equation*}
\Bar\nabla F(\ox):=\co\big\{\lim\nabla F(x_k)\;\big|\;x_k\to\ox,\;k\to\infty,\;F\;\mbox{ is differentiable at }\;x_k\big\}.
\end{equation*}\vspace*{-0.15in}

Now let $\vph\colon\R^n\to\oR$ be an extended-real-valued function with the domain and epigraph given by
\begin{equation*}
\dom\vph:=\big\{x\in\R^n\left|\right.\vph(x)<\infty\big\}\;\textrm{ and }\;\epi\vph:=\big\{(x,\al)\in\R^{n+1}\left|\right.\al\ge\vph(x)\big\},
\end{equation*}
respectively, where $\oR:=(-\infty,\infty]$. Assuming that $\ph$ is lower semicontinuous (l.s.c.) around $\ox\in\dom\ph$, its (first-order) {\em subdifferential} at $\ox$ is defined geometrically via the normal cone \eqref{nor_con} by
\begin{equation}\label{1sub}
\partial\vph(\ox):=\big\{v\in\R^n\left|\right.(v,-1)\in N\big((\ox,\vph(\ox));\epi\vph\big)\big\}.
\end{equation}
It is easy to see that in the case where $\ph(x):=\dd(x;\O)$ is the {\em indicator function} of a locally closed set $\O\subset\R^n$ that equals $0$ for $x\in\O$ and $\infty$ otherwise, we have the relationship $\partial\dd(\ox;\O)=N(\ox;\O)$ for all $\ox\in\O$, i.e., get back to the normal cone \eqref{nor_con}.\vspace*{-0.05in}

Observe finally that if $\ph\colon\R^n\to\oR$ is a locally Lipschitzian function around $\ox$, then the subdifferential $\partial\ph(\ox)$ is a nonempty compact set in $\R^n$. On the other hand, if $F\colon\R^n\to\R^m$ is a single-valued and locally Lipschitzian mapping around $\ox$, then we have the {\em coderivative scalarization} formula:
\begin{equation}\label{scal}
D^*F(\ox)(u)=\partial\la u,F\ra(\ox)\;\mbox{ for all }\;u\in\R^m.
\end{equation}\vspace*{-0.4in}

\section{Second-Order Constructions and Their Calculations}\label{2nd-order}
\setcounter{equation}{0}\vspace*{-0.1in}

Next we turn to the {\em second-order} generalized differential constructions for extended-real-valued functions that play a crucial role in our study of controlled sweeping processes with nonsmooth perturbations. This is due to the fact that the sweeping differential inclusion in \eqref{SP} are described by the normal cone mapping, while an adjoint system to it is naturally expressed via the coderivative \eqref{cod}, which is a generalized adjoint derivative operator. Such a {\em dual derivative-of-derivative} approach to second-order generalized differentiation was suggested in \cite{m92} and than has been strongly developed and applied in many publications; see, e.g., the books \cite{m-book1,m-book2,m18} with the references and commentaries therein.\vspace*{-0.05in}

Given $\ph\colon\R^n\to\oR$ with $\ox\in\dom\ph$, its {\em second-order subdifferential} (or {\em generalized Hessian}) at $\ox$ relative to $\ov\in\partial\ph(\ox)$ is defined as the mapping $\partial^2\ph(\ox,\ov)\colon\R^n\tto\R^n$ with the values
\begin{equation}\label{2sub}
\partial^2\ph(\ox,\ov)(u):=\big(D^*\partial\ph\big)(\ox,\ov)(u),\quad u\in\R^n,
\end{equation}
generated by the coderivative \eqref{cod} of the first-order subdifferential mapping $F:=\partial\ph$ from \eqref{1sub}. If $\ph:=N_\O$, we have $\partial^2\dd_\O(\ox,\ov)=D^*N_\O(\ox,\ov)$ for any $\ox\in\O$ and $\ov\in N_\O(\ox)$.

Now we consider the set-valued mapping $F\colon\R^n\times\R^d\tto\R^n$ that appears on the right-hand side of the sweeping differential inclusion \eqref{SP} as
\begin{equation}\label{F0}
F(x,u):=N(x;C)+g(x,u).
\end{equation}
It follows from the classical Motzkin theorem of the alternative that
\begin{equation}\label{F}
F(x,u):=\Big\{\sum_{j\in I(x)}\lm^j x^j_*\;\Big|\;\lm^j\ge 0\Big\}+g(x,u),
\end{equation}
where the collection of the active constraint indices $I(x)$ for the polyhedron \eqref{C} at $x\in C$ is taken from \eqref{aci}. Note that the coderivative of the mapping $F$ in \eqref{F0} relates to the second-order subdifferential of $\dd_C$, and that the coderivative of $g$ admits the subdifferential characterization \eqref{scal}. On the other hand, it is more convenient for us to work with {\em bounded} differential inclusions instead of the unbounded one \eqref{SP} in terms of the normal cone mapping. It is possible due to the remarkable result by Thibault from \cite[Theorem~3.1]{Thi}, which shows that the latter sweeping differential inclusion is equivalent to
\begin{equation}\label{bound}
-\dot{x}(t)\in N\big(x(t);C\big)\cap M\B+g\big(x(t),u(t)\big)\;\textrm{ a.e. }\;t\in[0,T]
\end{equation}
for any $M>0$ sufficiently large, where $\B$ stands for the unit ball in $\R^n$. Having in mind representation \eqref{F} of the mapping $F$ in terms of the generating vectors $x^j_*$ of the convex polyhedron \eqref{C}, consider the following subsets of the active constraint indices $I(\ox)$:
\begin{equation}\label{c56}
I_0(y):=\big\{j\in I(\ox)\;\big{|}\;\la x^j_*,y\ra=c_j\big\}\;\textrm{ and }\;I_>(y):=\big\{j\in I(\ox)\;\big{|}\;\la x^j_*,y\ra>c_j\big\},\quad y\in\R^n.
\end{equation}\vspace*{-0.2in}

Now we are ready to derive a crucial upper estimate of the coderivative \eqref{cod} of the sweeping control mapping \eqref{F}, which is a second-order subdifferential construction, entirely in terms of the given data of \eqref{SP}. This is done under our standing assumptions. If in addition the LICQ is imposed, we arrive at the precise equality formula, which is also used in some (while not major) results below. \vspace*{-0.1in}

\begin{theorem} {\bf(second-order calculation for nonsmooth sweeping processes).}\label{Thm6.1ns} Given $F$ from \eqref{F0} and \eqref{F}, fix $(\ox,\ou)\in\gph F$ and $M>0$ to be sufficiently large. Define the mappings $F_1\colon\R^n\tto\R^n$ and $F_2\colon\R^n\times\R^d\to\R^n$ by $F_1(x):=N(x;C)\cap M\B$ and $F_2(x,u):=g(x,u)$, respectively.
Then for any $(x,u)\in C\times U$ and $\o\in N(x;C)\cap(M\B)+g(x,u)$ we have the coderivative upper estimate
\begin{equation}\label{c57ns}
D^*F(x,u,\o)(w)\subset\bigg\{z\in\R^{n+d}\;\bigg|\;z\in\partial\la w,g\ra(x,u)+\bigg(\sum_{j\in I_0(w)\cup I_>(w)}\gg^j x^j_*,0\bigg)\bigg\},
\end{equation}
where $w\in\dom D^*N_C(x,\o-g(x,u))$, where $I_0(w)$ and $I_>(w)$ are taken from \eqref{c56} with $\ox=x$, and where $\gg^j\in\R$ for $j\in I_0(w)$ while $\gg^j\ge 0$ for $j\in I_>(w)$. If in addition the vectors $\{x^j_*\;|\;j\in I(\ox)\}$ are linearly independent, then the domain $\dom D^*F_1(x,\o-g(x,u))=\dom D^*N_C\big(x,\o-g(x,u)\big)$ can be computed by
\begin{equation}\label{cod-dom}
\dom D^*N_C\big(x,\o-g(x,u)\big)=\bigg\{w\bigg|\exists\lm^j\ge 0\;\textrm{ with }\;\o-g(x,u)=\underset{j\in I(x)}{\sum}\lm^j x^j_*,\;\lm^j>0\Rightarrow\la x^j_*,w\ra=0\bigg\}.
\end{equation}
If finally the mapping $g(x,u)$ is ${\cal C}^1$-smooth around $(\ox,\ou)$, the coderivative inclusion \eqref{c57ns} reduces to
\begin{equation}\label{c57smooth}
D^*F(x,u,\o)(w)\subset\bigg\{z=\bigg(\nabla_x g(x,u)^*w+\sum_{j\in I_0(w)\cup I_>(w)}\gg^j x^j_*,\nabla_u g(x,u)^*w\bigg)\bigg\},
\end{equation}
which holds as equality provided that the vectors $\{x^j_*\;|\;j\in I(\ox)\}$ are linearly independent.
\end{theorem}\vspace*{-0.1in}
{\bf Proof.} Let us first estimate the coderivative of the sum $F_1+F_2$ at $(x,u,\o)$ by using the sum rule from \cite[Theorem~3.9(ii)]{m18}. Observe that the sets
\begin{equation*}
S(x,u,\o):=\Big\{(y_1,y_2)\in\R^n\times\R^n\;\Big|\;y_1\in F_1(x),\;y_2\in F_2(x,u),\;y_1+y_2=\o\Big\}
\end{equation*}
in the aforementioned theorem reduce in our case to the form
\begin{equation}\label{S}
S(x,u,\o)=\Big\{\big(\o-g(x,u),g(x,u)\big)\in\R^n\times\R^n\;\Big|\;\o\in g(x,u)+N(x;C)\cap(M\B)\Big\}
\end{equation}
while being obviously uniformly bounded by the construction of $F_1$ and the assumptions on $U$ and $g$ in (H2) and (H3), respectively. The qualification condition in \cite[Theorem~3.9(ii)]{m18} reads as
\begin{equation*}
D^*F_1(x,y_1)(0)\cap\big(-D^*F_2(x,u,y_2)(0)\big)=\big\{(0,0)\big\},
\end{equation*}
and it holds by the assumed Lipschitz continuity of $g$ due to the necessity part of the coderivative criterion for the Lipschitz continuity from \cite[Theorem~3.3]{m18}, which ensures that $D^*F_2(x,u,y_2)(0)=\{0\}$ for all the triples $(x,u,y_2)$ under consideration. Applying now the coderivative sum rule from \cite[Theorem~3.9(ii)]{m18} to the sum $F_1+F_2$, we arrive at the inclusion
\begin{eqnarray}\label{sum1}
D^*(F_1+F_2)(x,u,\o)(w)\subset\underset{(y_1,y_2)\in S(x,u,\o)}{\bigcup}\bigg(D^*F_1(x,y_1)(w)+D^*F_2(x,u,y_2)(w)\bigg)
\end{eqnarray}
for all the corresponding quadruples $(x,u,\o,w)$. To proceed further, recall that $F_1=N_C\cap M\B$, where the constant $M>0$ can be chosen so large that $(x,y_1)\in{\rm int}(\gph M\B)=\R^n\times{\rm int}(M\B)$. This tells us that $N((x,y_1);\gph(M\B))=\{0,0\}$. Applying now the coderivative intersection rule from \cite[Proposition~3.20]{m-book1} to the mapping intersection $N_C\cap M\B$, we see that the qualification condition therein
\begin{equation*}
N\big((x,y_1);\gph N_C\big)\cap\gph(M\times\B)=\big\{(0,0)\big\}
\end{equation*}
holds automatically. The intersection formula from \cite[Proposition~3.20]{m-book1} reduces in our setting to
\begin{equation*}
D^*F_1(x,y_1)(w)=\Big\{D^*N_C(x,y_1)(w)\;\Big|\;w\in\dom D^*N_C\big(x,\o-g(x,u)\big)\Big\}
\end{equation*}
with $\dom D^*N_C(x,\o-g(x,u))=\dom D^*F_1(x,\o-g(x,u))$. Substituting the latter into the coderivative sum rule \eqref{sum1} with taking into account the form of the sets $S(x,u,\o)$ in \eqref{S}, the scalarization formula \eqref{scal} for the Lipschitzian mapping $g(x,u)$, and the second-order subdifferential definition \eqref{2sub}, we arrive at the following inclusion for the coderivative of the mapping $F$ from \eqref{F0}:
\begin{equation}\label{codF}
D^*F(x,u,\o)(w)\subset\bigg\{z\in\R^{n+d}\;\bigg|\;z\in\Big(\partial^2\dd_C\big(x,\o-g(x,u)\big)(w),0\Big)+\partial\la w,g\ra(x,u)\bigg\},
\end{equation}
where $\o\in N(x;C)\cap(M\B)+g(x,u)$, $(x,u)\in X\times U$, and $w\in\dom D^*N_C(x,\o-g(x,u))$. Now using the upper estimate of $\partial^2\dd_C$ in \eqref{codF} established in \cite[Theorem~4.5]{hmn}, we arrive at the main result \eqref{c57ns} of the theorem under the imposed standing assumptions. If furthermore  $g$ is ${\cal C}^1$-smooth around $(\ox,\ou)$, then the coderivative sum rule \eqref{sum1} holds as equality. In this case we have the decomposition
\begin{equation}\label{decomp}
\nabla g(x,u)=\nabla_x g(x,u)\times\nabla_u g(x,u)\;\mbox{ for all }\;(x,u)\;\mbox{ near }\;(\ox,\ou).
\end{equation}
Combining \eqref{decomp} with the aforementioned estimate of $\partial^2\dd_C$ from \cite[Theorem~4.5]{hmn} verifies the inclusion in \eqref{c57smooth}. Since the latter estimate of $\partial^2\dd_C$ is proved in \cite[Theorem~4.5]{hmn} to become an equality under the additional LICQ assumption, we deduce from the above that the equality holds in \eqref{c57smooth} under the LICQ and smoothness of $g$. The domain formula \eqref{cod-dom} under LICQ follows from \cite[Theorem~4.5]{hmn}. $\h$\vspace*{-0.03in}

Note that the coderivative and subdifferential counterparts of the decomposition formula \eqref{decomp} hold, as the desired inclusions ``$\subset$", under some restrictive regularity assumptions; see \cite[Corollaries~3.17 and 3.44]{m-book1}. We prefer to avoid them and thus employ in this paper the upper estimate \eqref{c57ns}, which is sufficient for deriving the main necessary optimality conditions for the discrete-time and continuous-time problems.\vspace*{-0.2in}

\section{Necessary Conditions for Discrete Approximations}\label{disc}
\setcounter{equation}{0}\vspace*{-0.1in}

The method of discrete approximations in the derivation of necessary conditions for optimal solutions to continuous-time problems consists of constructing well-posed discrete-time problems that adequately approximate a given local minimizer of the original one, establishing necessary optimality conditions for approximating discrete optimal solutions, and then justifying the passage to the limit from the discrete optimality conditions to derive the desired ones for the local minimizer in question. This approach was initiated in \cite{m95} for Lipschitzian differential inclusions (see also the book \cite{m-book2} with the commentaries therein) and then strongly developed for various versions of controlled sweeping processes in \cite{cm1,cm2,cm3,chhm1,chhm3,cmn1,cmn,hm} and other publications. In this paper we follow the scheme of \cite{cmn1,cmn} in the construction of discrete approximations and establishing the appropriate convergence of discrete optimal solutions. Then we derive novel necessary optimality conditions for discrete-time systems by using the tools of generalized differentiation overviewed in Section~\ref{tools} and particularly the new second-order results of Theorem~\ref{Thm6.1ns}.\vspace*{-0.03in}

Consider the standard {\em Euler explicit scheme} for the replacement of the time derivative in \eqref{Problem} by
\begin{equation*}
\dot x(t)\approx\frac{x(t+h)-x(t)}{h}\textrm{ as }h\dn 0
\end{equation*}
and formalize it as follows. For any $m\in\N:=\{1,2,\ldots\}$ denote the discrete mesh on $[0,T]$ by
\begin{equation*}
\DD_m:=\big\{0=t^0_m<t^1_m<\ldots<t^{2^m}_m=T\big\}\textrm{ with }h_m:=t^{i+1}_m-t^i_m=\frac{T}{2^m},
\end{equation*}
which is made uniform for simplicity, and define the sequence of {\em discrete-time controlled sweeping processes}
\begin{equation}\label{e:3.4}
-x^{i+1}_m\in-x^i_m+h_m F(x_m^{i},u_m^{i}),\quad i=0,\ldots,2^m-1,\quad x^0_m:=x_0\in C,
\end{equation}
with {\em control constraints} $u_m^{i}\in U$ for all $i=0,\ldots,2^m-1$, where $F$ is taken from \eqref{F0}. The latter implies that we automatically have the {\em state constraints} at each discrete time:
\begin{equation}\label{disc-state}
x^i_m\in C\Longleftrightarrow\la x^j_*,x^i_m\ra\le c_j\;\mbox{ whenever }\;j=1,\ldots,s,\quad i=0,\ldots,2^m-1,\quad m\in\N.
\end{equation}\vspace*{-0.2in}

Let $(\ox(\cdot),\ou(\cdot))$ be {\em any feasible} pair for the original problem $(P)$ such that $\ou(\cdot)$ is a function of {\em bounded variation} (BV) on $[0,T]$ admitting a {\em right continuous representative}, i.e., a function with these properties that coincides with $\ou(\cdot)$ for a.e.\ $t\in[0,T]$. In what follows we use the same notation for such a representative and add the aforementioned properties to the {\em standing assumptions} on the local optimal control $\ou(\cdot)$ to $(P)$ under consideration. Denote $I_m^i:=[t_m^{i-1},t_m^i)$ for $i=0,\ldots,2^m-1$.\vspace*{-0.03in}

It is proved in \cite[Theorem~3.1]{cmn} that, given such a feasible solution $(\ox(\cdot),\ou(\cdot))$ to $(P)$, for each $m\in\mathbb{N}$ there exist state-control pairs $(x_m(t),u_m(t))$
and perturbation terms $r_m(t)\ge 0$ and $\rho_m(t)\in\mathbb{B}$ defined on $[0,T]$ and having the following properties:
\begin{equation*}
\dot{x}_m(t)\in-N\big(x_m(t_m^{i});C\big)-g\big(x_m(t_m^{i}),u_m(t)\big)+r_m(t)\rho_m(t),
\end{equation*}
\begin{equation*}
x_m(t_m^i)=\ox(t_m^i)\in C\;\textrm{ for each }\;i=1,\ldots,2^m\;\textrm{ with }\;x_m(0)=x_0
\end{equation*}
for all $t\in(t_m^{i-1},t_m^i)$ and $i=1,\ldots,2^m$, where $x_m\colon[0,T]\to\R^n$ are continuous on $[0,T]$, affine on each interval $I_m^i$, and converge strongly in
$W^{1,2}([0,T];\R^n)$ to $\ox(\cdot)$; where control mappings $u_m\colon[0,T]\to U$ are constant on $I_m^i$ and converge to $\ou(\cdot)$ strongly in $L^2([0,T];\R^d)$ and pointwise on $[0,T]$; and where $r_m\colon[0,T]\to[0,\infty)$ and $\rho_m\colon[0,T]\to\mathbb{B}$ are constant on each interval $I_m^i$ with $r_m(\cdot)\to 0$ in $L^2([0,T];\R)$.

From now on, the above feasible pair $(\ox(\cdot),\ou(\cdot))$ is a given {\em relaxed $W^{1,2}\times L^2$-local minimizer} of $(P)$. Employing the aforementioned approximation result for this feasible pair, let us construct a sequence of {\em discrete-time sweeping optimal control} problems $(P_m)$ that provides a desired {\em strong approximation} of $(\ox(\cdot),\ou(\cdot))$ by {\em optimal solutions} to $(P_m)$. Given $\ve>0$ in Definition~\ref{relaxed} of relaxed $W^{1,2}\times L^2$-local minimizers, for each $m\in\N$ we define $(P_m)$ as follows:
\begin{equation*}
\textrm{minimize }\;J_m[x_m,u_m]:=\vph\big(x_m(T)\big)+\frac{1}{2}\sum_{i=0}^{2^m-1}\int_{t^i_m}^{t^{i+1}_m}\(\n\frac{x^{i+1}_m-x^i_m}{h_m}-\dot{\ox}(t)\en^2+
\n u^{i}_m-\ou(t)\en^2\)dt
\end{equation*}
over discrete functions $(x_m,u_m):=(x^0_m,x^1_m,\ldots,x^{2^m}_m,u^0_m,u^1_m,\ldots,u^{2^m-1}_m)$ satisfying the constraints
\begin{equation}\label{SP-disc}
x^i_m-x^{i+1}_m\in h_m F_m(t^i_m,x_m^i,u_m^i)\;\textrm{ for }\;i=0,\ldots,2^m-1,
\end{equation}
where $F_m(t,x,u):=F(x,u)-r_m(t)\rho_m(t)$, and where
\begin{equation*}
\langle x_\ast^{j},x^{2^m}_m\rangle\le c_{j}\;\textrm{ for all }\;j=1,\ldots,s\;\text{ with }\;x^0_m:=x_0\in C,
\;u^0_m:=\ou (0),
\end{equation*}
\begin{equation*}
\sum_{i=0}^{2^m-1}\int_{t^i_m}^{t^{i+1}_m}\(\n\frac{x^{i+1}_m-x^i_m}{h_m}-\dot{\ox}(t)\en^2+\n u_m^i-
\ou(t)\en^2\)dt\le\frac{\e}{2},
\end{equation*}
\begin{equation*}
u^i_m\in U\;\textrm{ for }\;i=0,\ldots,2^m-1,
\end{equation*}
and where the state constraints \eqref{disc-state} are implicitly included in \eqref{SP-disc} due to the structure of $F$.

The following theorem is taken from \cite[Proposition~4.1 and Theorem~4.2]{cmn}.\vspace*{-0.1in}

\begin{theorem}{\bf(strong convergence of discrete optimal solutions).}\label{ThmStrong} Let $(\ox(\cdot),\ou(\cdot))$ be a relaxed $W^{1,2}\times L^2$-local minimizer for the sweeping control problem $(P)$. Then each problem $(P_m)$ has an optimal solution whenever $m\in\N$ is sufficiently large. Furthermore, any sequence of optimal solutions $(\ox_m(\cdot),\ou_m(\cdot))$ to $(P_m)$, which is extended to $[0,T]$ piecewise linearly for $\ox_m(\cdot)$ and piecewise constantly for $\ou_m(\cdot)$, converges strongly to $(\ox(\cdot),\ou(\cdot))$ as $m\to\infty$ in the norm topology of $W^{1,2}([0,T];\R^n)\times L^2([0,T];\R^d)$.
\end{theorem}\vspace*{-0.07in}

Our next major step is to obtain necessary conditions for optimal solutions to each problem $(P_m)$. To proceed, we employ the nonconvex generalized differentiation tools of variational analysis discussed in Sections~\ref{tools} and \ref{2nd-order} with the main impact of the new second-order calculations of Theorem~\ref{Thm6.1ns}. \vspace*{-0.03in}

In what follows we derive two results in this direction that both use Theorem~\ref{Thm6.1ns}. The first result provides necessary optimality conditions for $(P_m)$ that involve the normal cone to the graph of the velocity mapping $F$ from \eqref{F}. The obtained relationships can be treated as discrete counterparts of the Euler-Lagrange conditions in our setting. Theorem~\ref{Thm6.1ns} is very instrumental in the proof of these conditions to justify the application of the basic normal cone calculus from \cite{m18}. The second (main) result here establishes the collection of necessary optimality conditions for $(P_m)$ expressed entirely via the given problem data. It becomes possible due to the explicit second-order calculations of Theorem~\ref{Thm6.1ns} and leads us in Section~\ref{optim} to deriving necessary optimality conditions for the relaxed $W^{1,2}\times L^2$-local minimizer $(\ox(\cdot),\ou(\cdot))$ to the original problem $(P)$ by passing to the limit from discrete approximations. \vspace*{-0.1in}

\begin{theorem}{\bf(necessary conditions of the extended Euler-Lagrange type for discrete sweeping processes).}\label{Thm5.2*} Let $(\ox_m,\ou_m)=(\ox^0_m,\ldots,\ox^{2^m}_m,\ou^0_m,\ldots,\ou^{2^m-1}_m)$ be an optimal solution to problem $(P_m)$ for any fixed $m\in\N$, where $F$ is defined in \eqref{F0} under the fulfillments of the standing assumptions around $\ox:=\ox_m$. Then there are dual elements $\lm_m\ge 0,\;\psi_m=(\psi^0_m,\ldots,\psi^{2^m-1}_m)\in\R^{2^m}_+$ with $\psi^i_m\in N\(\ou^i_m;U\)$ for $i=0,\ldots,2^m-1$, as well as $\xi_m=(\xi^1_m,\ldots,\xi^s_m)\in\R^s_+$ and $p^i_m\in\R^n$ for $i=0,\ldots,2^m$ such that
\begin{equation}\label{e:5.8*}
\lm_m+\n\xi_m\en+\sum_{i=0}^{2^m-1}\n p^{i}_m\en+\n\psi_m\en\ne 0,
\end{equation}
\begin{equation}\label{xi1}
\xi^j_m\(\la x^j_*,x^{2^m}_m\ra-c_j\)=0,\quad j=1,\ldots,s,
\end{equation}
\begin{equation} \label{mutx}
-p^{2^m}_m=\lm_m\vth^{2^m}_m+\sum_{j=1}^s\xi^{j}_m x^{j}_\ast,\;\textrm{ with }\;\vth^{2^m}_m\in\partial\vph(\ox^{2^m}_m),
\end{equation}
\begin{equation}\label{e:5.10*}
\begin{array}{ll}
&\disp\bigg(\frac{p^{i+1}_m-p^{i}_m}{h_m},-\frac{1}{h_m}\lm_m\th^{iu}_m,\frac{1}{h_m}\lm_m\th^{iy}_m-p^{i+1}_m\bigg)\\\\
&\disp\in\bigg(0,\dfrac{1}{h_m}\psi^i_m,0\bigg)+N\bigg(\bigg(\ox^i_m,\ou^i_m,-\frac{\ox^{i+1}_m-\ox^i_m}{h_m}\bigg);\;\gph F_m\bigg),\quad i=0,\ldots,2^m-1,
\end{array}
\end{equation}
with the auxiliary vectors $\th^{iy}_m$ and $\th^{iu}_m$ in \eqref{e:5.10*} defined by
\begin{equation}\label{theta}
\(\th^{iy}_m,\th^{iu}_m\):=\bigg(\int_{t^i_{m}}^{t^{i+1}_{m}}\bigg(\frac{\ox^{i+1}_m-\ox^i_m}{h_m}-\dot{\ox}(t)\bigg)dt,\int_{t^i_{m}}^{t^{i+1}_{m}}\(\ou^i_m-\ou(t)\)dt\bigg).
\end{equation}
\end{theorem}\vspace*{-0.05in}
{\bf Proof.} Fix the number $\e>0$ from Definition~\ref{relaxed} of the given $W^{1,2}\times L^2$-local minimizer $(\ox(\cdot),\ou(\cdot))$ to the original problem $(P)$ and consider the ``long" vector
\begin{equation*}
z:=(x^0_m,\ldots,x^{2^m}_m,u^0_m,\ldots,u^{2^m-1}_m,y^0_m,\ldots,y^{2^m-1}_m)\in\R^{(2\cdot 2^m+1)n+2^m\cdot d}
\end{equation*}
with the fixed starting point $x^0_m$. It is easy to see that the discrete-time problem $(P_m)$ for each fixed $m\in\N$ can be written as the following equivalent problem of nondynamic problem of {\em mathematical programming} $(MP)$ with respect to the variable $z$:
\begin{equation*}
\textrm{minimize }\;\phi_0(z):=\vph\big(x(T)\big)+\frac{1}{2}\sum_{i=0}^{2^m-1}\int_{t^i_m}^{t^{i+1}_m}\n\(y^i_m-
\dot{\ox}(t),u_m^i-\ou(t)\)\en^2dt
\end{equation*}
subject to the smooth equality and inequality constraints
\begin{equation}\label{phi}
\phi(z):=\sum_{i=0}^{2^m-1}\int_{t^i_m}^{t^{i+1}_m}\n\(y^i_m,u^i_m\)-\big(\dot{\ox}(t),\ou(t)\big)\en^2dt-\frac{\e}{2}\le 0,
\end{equation}
\begin{equation}\label{g}
g_i(z):=x^{i+1}_m-x^i_m-h_m y^i_m=0,\quad i=0,\ldots,2^m-1,
\end{equation}
\begin{equation}\label{h}
h_{j}(z):=\la x^{j}_\ast,x^{2^m}_m\ra-c_{j}\le 0,\quad j=1,\ldots,s,
\end{equation}
as well as the increasingly many geometric constraints
\begin{equation}\label{gra}
z\in\Xi_i:=\left\{(x^0_m,\ldots,y^{2^m-1}_m)\in\R^{(2\cdot2^m+1)n+2^m\cdot d}\;\Big|-y^i_m
\in F_m(t^i_m,x^i_m,u^i_m)\right\},\;i=0,\ldots,2^m-1,
\end{equation}
\begin{equation}\label{xi}
z\in\Xi_{2^m}:=\Big\{(x^0_m,\ldots,y^{2^m-1}_m)\in\R^{(2\cdot 2^m+1)n+2^m\cdot d}\;\Big|\;x^0_m\;\textrm{ is fixed}\Big\},
\end{equation}
\begin{equation}\label{om}
z\in\O_i:=\Big\{(x^0_m,\ldots,y^{2^m-1}_m)\in\R^{(2\cdot 2^m+1)n+2^m\cdot d}\;\Big|\;u^i_m\in U\Big\},\quad i=0,\ldots,2^m-1.
\end{equation}
Nonstandard constraints in $(MP)$ are geometric ones among which the {\em graphical} constraints \eqref{gra} are the most challenging. Luckily, we can be handled them in terms of our basic normal cone \eqref{nor_con} (i.e., via the coderivative of $F_m$), while definitely not via its convexification.\vspace*{-0.03in}

To proceed, let us apply the appropriate necessary optimality conditions developed in nonsmooth constrained optimization to the optimal solution $\oz:=\(\ox^0_m,\ldots,\ox^{2^m}_m,\ou^0_m,\ldots,\ou^{2^m-1}_m,\oy^0_m,\ldots,\oy^{2^m-1}_m\)$ for $(MP)$ corresponding to the optimal solution $(\ox_m,\ou_m)$ for the discrete-time control problem $(P_m)$. Using \cite[Theorem~6.5]{m18} and supporting it by the intersection rule for the normal cone \eqref{nor_con} taken from \cite[Theorem~2.16 and Corollary~6.17]{m18} with including the qualification conditions therein into the nontriviality relation, we find $\lm_m\ge 0$, $\xi_m=(\xi^{1}_m,\ldots,\xi^{s}_m)\in\R^s_+$, $p^i_{m}\in\R^{n}$ as $i=1,\ldots,2^m$, and
\begin{equation*}
z^*_i=\big(x^*_{0i},\ldots,x^*_{2^mi},u^*_{0i},\ldots,u^*_{(2^m-1)i},
y^*_{0i},y^*_{1i},\ldots,y^*_{(2^m-1)i}\big),\quad i=0,\ldots,2^m,
\end{equation*}
which are not equal to zero simultaneously and satisfy the conditions
\begin{equation}\label{69}
z^*_i\in\left\{\begin{matrix}
N(\oz;\Xi_i)+N(\oz;\O_i)\;\textrm{ if }\;i\in\big\{0,\ldots,2^m-1\big\},\\
N(\oz;\Xi_i)\;\textrm{ if }\;i=2^m,
\end{matrix}\right.
\end{equation}
\begin{equation}\label{disc-lagr}
-z^*_0-\ldots-z^*_{2^m}\in\lm_m\partial\phi_0(\oz)+\sum_{j=1}^{s} \xi^{j}_m\nabla h_{j}(\oz)+
\sum_{i=0}^{2^m-1}\nabla g_i(\oz)^*p^{i+1}_m,
\end{equation}
\begin{equation}\label{71+}
\xi^{j}_m h_{j}(\oz)=0,\quad j=1,\ldots,s,
\end{equation}
for all $m\in\N$ sufficiently large. Note that the Lagrangian condition \eqref{disc-lagr} does not contain the term reflecting the inequality constraint \eqref{phi}. It is due to Theorem~\ref{ThmStrong} above ensuring that this constraint is inactive for sufficient large $m$, and thus the corresponding multiplier is zero. Observe also that the upper expression in \eqref{69} benefits from the application of the basic intersection rule for basic normals taken from \cite[Theorem~2.16]{m18}, which tells us that
\begin{equation*}
N(\oz;\Xi_i\cap\O_i)\subset N(\oz;\Xi_i)+N(\oz;\O_i)\;\mbox{ whenever }\;i=0,\ldots,2^m-1
\end{equation*}
under the validity of the normal qualification condition
\begin{equation}\label{qc}
N(\oz;\Xi_j)\cap\big(-N(\oz;\O_i)\big)=\{0\}
\end{equation}
for each $i=0,\ldots,2^m-1$. To check the fulfillment of \eqref{qc}, fix any $i\in\{0,\ldots,2^m-1\}$ and pick $z^*_i\in N(\oz;\Xi_j)\cap(-N(\oz;\O_i))$. Then we have by the structures of $\Xi_i$ in \eqref{gra} and $\O_i$ in \eqref{om} that
\begin{equation}\label{qc0}
\big(x^*_{ii},u^*_{ii},-y^*_{ii}\big)\in N\Big(\Big(\ox^i_m,\ou^i_m,-\dfrac{\ox^{i+1}_m-\ox^i_m}{h_m}\Big);\gph F_m\Big)\;\mbox{ and }\;-u^*_{ii}\in N\(\ou^i_m;U\)
\end{equation}
with $\gph F_m:=\gph F_m(t^i_m\cdot,\cdot)$, while all the other components of $z^*_i$ are zero. It is obvious from \eqref{qc0} that $x^*_{ii}=0$ and $y^*_{ii}=0$ for these components of $z^*_i$. Thus \eqref{qc} reduces to the implication
\begin{equation}\label{qc1}
\big(0,u^*_{ii},0\big)\in N\Big(\Big(\ox^i_m,\ou^i_m,-\dfrac{\ox^{i+1}_m-\ox^i_m}{h_m}\Big);\gph F_m\Big)\Longrightarrow u^*_{ii}=0,\quad j=0,\ldots,2^m-1.
\end{equation}
Using the relationship between $F_m$ and $F$ in \eqref{SP-disc} and the coderivative definition \eqref{cod}, we rewrite the left-hand side of the implication in \eqref{qc1} as
\begin{equation}\label{qc2}
(0,u^*_{ii})\in D^*F\Big(\ox^i_m,\ou^i_m,r_m(t^i_m)\rho_m(t^i_m)-\frac{x^{i+1}_m-\ox^i_m}{h_m}\Big)\big(0\big).
\end{equation}
By employing the coderivative estimate \eqref{c57ns} of Theorem~\ref{Thm6.1ns}, it follows from \eqref{qc2} that
\begin{equation*}
(0,u^*_{ii})\in\partial\la 0,g\ra(\ox^i_m,\ou^i_m)+\bigg(\sum_{j\in I_0(0)\cup I_>(0)}\gg^j x^j_*,0\bigg),
\end{equation*}
which obviously yields $u^*_{ii}=0$ and thus verifies the qualification condition \eqref{qc}.\vspace*{-0.03in}

Invoking again the structures of $\Xi_i$ and $\O_i$, we rewrite the upper formula in \eqref{69} as
\begin{equation}\label{e:5.18*}
\big(x^*_{ii},u^*_{ii}-\psi^{i}_m,-y^*_{ii}\big)\in N\Big(\Big(\ox^i_m,\ou^i_m,-\frac{\ox^{i+1}_m-\ox^i_m}{h_m}\Big);\gph F_m\Big)\;\textrm{ for }\;i=0,\ldots,2^m-1
\end{equation}
with $\psi^{i}_m\in N(\ou^i_m;U)$ for all $i=0,\ldots,2^m-1$ and with every other components of $z^*_i$ being zero. It easily follows from the structure of $\Xi_{2^m}$ in \eqref{xi} that $x^*_{02^m}$ and $u^*_{02^m}$ are the only nonzero components of the vector $z^*_{2^m}$ in the lower formula of \eqref{69}. Thus we deduce from
\eqref{69} and \eqref{71+} that\vspace*{-0.05in}
\begin{equation}\label{inc1}
\begin{array}{ll}
\disp-z^*_0-\ldots-z^*_{2^m}\in\lm_m\partial\phi_0(\oz)+\sum_{j=1}^{s}\xi^{j}_m\nabla h_{j}(\oz)+\sum_{i=0}^{2^m-1}\nabla g_i(\oz)^*p^{i+1}_m\\
\disp\mbox{with }\;\xi^{j}_m\big(\la z^{j2^m}_m, x^{2^m}_m\ra-c^{j2^m}_m\big)=0\;\mbox{ for all }\;j=1,\ldots,s,
\end{array}
\end{equation}
where the latter follows from \eqref{71+}. The smoothness of $g_i$ in \eqref{g} and $h_j$ in \eqref{h} yields
\begin{equation*}
\Big(\sum_{i=0}^{2^m-1}\nabla g_i(\oz)^*p^{i+1}_m\Big)_{x^i_m}=\left\{\begin{matrix}
-p^{1}_m&\mbox{ if }\;i=0,\\
p^{i}_m-p^{i+1}_m&\mbox{ if }\;i=1,\ldots,2^m-1,\\
p^{2^m}_m&\mbox{ if }\;i=2^m,
\end{matrix}\right.
\end{equation*}
\begin{equation*}
\Big(\sum_{i=0}^{2^m-1}\nabla g_i(\oz)^*p^{i+1}_m\Big)_{y^i_m}=\(-h_m p^{1}_m,-h_m p^{2}_m,\ldots,-h_m p^{2^m}_m\),
\end{equation*}
\begin{equation*}
\Big(\sum_{j=1}^{s}\xi^{j}_m\nabla h_{j}(\oz)\Big)_{x^{2^m}_m}=\Big(\sum_{j=1}^{s}\xi^{j}_m x_*^j\Big).
\end{equation*}
Furthermore, using the smoothness of the second term in \eqref{phi} and the calculation of its gradients together with the elementary subdifferential sum rule from \cite[Proposition~1.30(ii)]{m18} implies that
\begin{equation*}
\partial\phi_0(\oz)=\partial\vph(\ox^m_m)+\frac{1}{2}\sum_{i=0}^{2^m-1}\nabla\rho_i(\oz)\;\textrm{ with }\;
\rho_i(\oz):=\int_{t^i_m}^{t^{i+1}_m}\Big\|\Big(\frac{\ox^{i+1}_m-\ox^i_m}{h_m}-\dot{\ox}(t),\ou^i_m-\ou(t)\Big)\Big\|^2dt.
\end{equation*}
This allows us to express the subgradient set $\lm_m\partial\phi_0(\oz)$ in the form
\begin{equation*}
\lm_m\big(0,\ldots,0,\vth^{2^m}_m,\th^{0u}_m,\ldots,\th^{(2^m-1)u}_m,\th^{0y}_m,\ldots,\th^{(2^m-1)y}_{m}\big)\;
\textrm{ with }\;\vth^{2^m}_m\in\partial\vph(\ox^{2^m}_m)
\end{equation*}
\begin{equation*}
\mbox{with }\;(\th^{iu}_m,\th^{iy}_m)=\Big(\int_{t^i_m}^{t^{i+1}_m}\Big(\ou^i_m -\ou(t)\Big)dt,\;
\int_{t^i_m}^{t^{i+1}_m}\Big(\frac{\ox^{i+1}_m-\ox^i_m}{h_m}-\dot{\ox}(t)\Big)dt\Big),\quad i=0,\ldots,2^m-1.
\end{equation*}
Substituting the above calculations into the inclusion of \eqref{inc1}, we arrive at the relationships
\begin{equation}\label{daux}
-x^*_{00}-x^*_{02^m}=-p^{1}_m,
\end{equation}
\begin{equation}\label{e:5.21*}
-x^*_{ii}=p^{i}_m-p^{i+1}_m,\quad i=1,\ldots,2^m-1,
\end{equation}
\begin{equation}\label{e:5.24*}
0=\lm_m\vth^{2^m}_m+p^{2^m}_m+\sum_{j=1}^{s}\xi^{j}_m x_*^j\;\textrm{ with }\;
\vth^{2^m}_m\in\partial\vph(\ox_m^{2^m}),
\end{equation}
\begin{equation}\label{e:5.22*}
-u^*_{00}=\lm_m\th^{0u}_m\;\textrm{ and }\;-u^*_{ii}=\lm_m\th^{iu}_m,\quad i=1,\ldots,2^m-1,
\end{equation}
\begin{equation}\label{e:5.23*}
-y^*_{ii}=\lm_m\th^{iy}_m-h_m p^{i+1}_m,\quad i=0,\ldots,2^m-1.
\end{equation}\vspace*{-0.2in}

Now we are ready to complete the proof of the theorem. Note first that the complementary slackness conditions in \eqref{xi1} follow from the equalities in \eqref{inc1} and that the transversality condition \eqref{mutx} follows from $\eqref{e:5.24*}$. Then we get from \eqref{e:5.21*}, \eqref{e:5.22*}, and \eqref{e:5.23*} that
\begin{equation*}
\frac{x^*_{ii}}{h_m}=\frac{p^{i+1}_m-p^{i}_m}{h_m},\;\;\frac{u^*_{ii}}{h_m}=-\frac{1}{h_m}\lm_m\th^{iu}_m,\;\textrm{ and }\;\frac{y^*_{ii}}{h_m}
=-\frac{1}{h_m}\lm_m\th^{iy}_m+p^{i+1}_m.
\end{equation*}
Plugging the latter into the left-hand side of \eqref{e:5.18*} and defining $p^0_m:=x^*_{02^m}$ in addition to $p^i_m$ for $i=1,\ldots,2^m$, we come up to the discrete Euler-Lagrange inclusion \eqref{e:5.10*} with $(\th^{iy}_m,\th^{iu}_m)$ taken from \eqref{theta}.\vspace*{-0.03in}

It remains to verify the nontriviality condition \eqref{e:5.8*}. Suppose on the contrary that $\lm_m=0,\xi_m=0,\psi_m=0$, and $p^{i}_m=0$ as $i=0,\ldots,2^m-1$ yielding $x^*_{02^m}=p^{0}_{m}=0$. It follows from \eqref{e:5.24*} that $p^{2^m}_m=0$, which tells us that $p^{i}_m=0$ whenever $i=0,\ldots,2^m$. Combining \eqref{daux} and \eqref{e:5.21*} implies that $x^*_{ii}=0$ for all $i=0,\ldots,2^m-1$. Furthermore, we deduce from \eqref{e:5.22*} that $u^*_{ii}=0$ for all $i=0,\ldots,2^m-1$. Observe also from \eqref{e:5.23*} that $y^*_{ii}=0$ for all
$i=0,\ldots,2^m-1$. Taking into account that all the components of $z^*_i$ but $(x^*_{ii},u^*_{ii},y^*_{ii})$ are zero whenever $i=0,\ldots,2^m-1$, we conclude that $z^*_{i}=0$ for $i=0,\ldots,2^m-1$ and that similarly $z^*_{2^m}=0$. Thus $z^*_i=0$ for all $i=0,\ldots,2^m$. This contradicts the nontriviality condition for $(MP)$ and hence completes the proof of the theorem. $\h$\vspace*{-0.03in}

The final and most important result of this section provides necessary optimality conditions for any optimal solution to $(P_m)$ expressed entirely in terms of the problem data. It is based on the general Euler-Lagrange optimality conditions obtained in Theorem~\ref{Thm5.2*} via the normal cone to the graph of $F_m$ from \eqref{SP-disc} and on the second-order calculations of Theorem~\ref{Thm6.1ns} for the original mapping $F$ defined in \eqref{F0} with its equivalent representation in \eqref{F}. The obtained optimality conditions for the discrete problems $(P_m)$ are of their own interest, while being a vehicle for deriving the constructive necessary conditions for relaxed local minimizers of the sweeping control problem $(P)$ established in Section~\ref{optim}.\vspace*{-0.07in}

\begin{theorem}{\bf(necessary optimality conditions for discrete sweeping processes via initial data).}\label{Nonsmooth} Let $(\ox_m,\ou_m)$ be an optimal solution to problem $(P_m)$ for each fixed $m\in\N$. Then there exist $\lm_m\ge 0$, $\psi^i_m\in N(\ou^i_m;U)$ as $i=0,\ldots,2^m-1$, and $p^i_m\in\R^n$ as $i=0,\ldots,2^m$ together with vectors $\eta^i_m\in\R^s_+$ as $i=0,\ldots,2^m$ and $\gg^i_m\in\R^s$ as $i=0,\ldots,2^m-1$ satisfying the following conditions:\\[1ex]
$\bullet$ {\sc Primal-dual dynamic relationships:}
\begin{equation}\label{ns87}
r_m(t^i_m)\rho_m(t^i_m)-\frac{\ox^{i+1}_m-\ox^i_m}{h_m}-g(\ox^i_m,\ou^i_m)=\sum_{j\in I(\ox^i_m)}\eta^{ij}_mx^j_*,
\end{equation}
\begin{equation}\label{nsconx}
\begin{array}{ll}
\disp\Big(\frac{p^{i+1}_m-p^{i}_m}{h_m},-\frac{1}{h_m}\lm_m\th^{iu}_m-\frac{1}{h_m}\psi^i_m\Big)\in\partial\Big\la-\frac{1}{h_m}\lm_m\th^{iy}_m+p^{i+1}_m,g\Big\ra (\ox^i_m,\ou^i_m)\\[1ex]
\disp+\bigg(\sum_{j\in I_0\(p^{i+1}_m-\frac{1}{h_m}\lm_m\th^{iy}_m\)\cup I_>\(p^{i+1}_m-\frac{1}{h_m}\lm_m\th^{iy}_m\)}\gg^{ij}_mx^j_*,\;0\bigg)\;\mbox{for all}\;i=0,\ldots,2^m-1\;\mbox{and}\; j=1,\ldots,s.
\end{array}
\end{equation}
$\bullet$ {\sc Transversality relationships:}
\begin{equation} \label{nsnmutx}
-p^{2^m}_m=\lm_m\vth^{2^m}_m+\sum_{j=1}^s\eta^{2^mj}_m x^j_\ast\;\textrm{ with }\;\vth^{2^m}_m\in\partial\vph(\ox^{2^m}_m).
\end{equation}
$\bullet$ {\sc Complementarity conditions:}
\begin{equation}\label{nseta}
\Big[\la x^j_\ast,\ox^i_m\ra<c_j\Big]\Longrightarrow\eta^{ij}_m=0\;\mbox{ for all }\;i=0,\ldots,2^m\;\mbox{ and }\;j=1,\ldots,s,
\end{equation}
\begin{equation}\label{ns94}
\[\la x^j_*,\ox^i_m\ra<c_j\]\Longrightarrow\gg^{ij}_m=0\;\textrm{ for }\;i=0,\ldots,2^m-1\;\textrm{ and }\;j=1,\ldots,s,
\end{equation}
\begin{equation}\label{ns93}
\left\{\begin{matrix}
j\in I_>\(p^{i+1}_m-\frac{1}{h_m}\lm_m\th^{iy}_m\)\Longrightarrow\gg^{ij}_m\ge 0,\\
\[j\notin I_0\(p^{i+1}_m-\frac{1}{h_m}\lm_m\th^{iy}_m\)\cup I_>\(p^{i+1}_m-\frac{1}{h_m}\lm_m\th^{iy}_m\)\]\Longrightarrow\gg^{ij}_m=0
\end{matrix}\right.
\end{equation}
for $i=0,\ldots,2^m-1$ and $j=1,\ldots,s$. If the vectors $\{x^j_*\;|\;j\in I(\ox^i_m)\}$ are linearly independent, then
\begin{equation}\label{ns96}
\eta^{ij}_m>0\Longrightarrow\Big[\Big\la x^j_*,p^{i+1}_m-\disp\frac{1}{h_m}\lm_m\th^{iy}_m\Big\ra=0\Big]\;\mbox{ as }\;i=0,\ldots,2^m-1,\quad j=1,\ldots,s.
\end{equation}
$\bullet$ {\sc Nontriviality condition:} We always have
\begin{equation}\label{nsentc}
\lm_m+\n p^{2^m}_m\en+\n\psi_m\en+\n\gg_m\en\ne 0.
\end{equation}
\end{theorem}\vspace{-0.05in}
{\bf Proof.} Note first that the aforementioned conditions on $\lm_m$ and $\psi^i_m$ come directly from Theorem~\ref{Thm5.2*}. Combining further the discrete Euler-Lagrange inclusion \eqref{e:5.10*} of Theorem~\ref{Thm5.2*} with the coderivative definition \eqref{cod}, we rewrite \eqref{e:5.10*} in the coderivative form
\begin{equation*}
\(\frac{p^{i+1}_m-p^{i}_m}{h_m},-\frac{1}{h_m}\lm_m\th^{iu}_m-\frac{1}{h_m}\psi^i_m\)\in D^*F_m\(\ox^i_m,\ou^i_m,-\frac{\ox^{i+1}_m-\ox^i_m}{h_m}\)\(-\frac{1}{h_m}\lm_m\th^{iy}_m+p^{i+1}_m\)
\end{equation*}
for all $i=0,\ldots,2^m-1$. Using the construction of the mapping $F_m$ in \eqref{e:5.10*} gives us
\begin{equation*}
r_m(t^i_m)\rho_m(t^i_m)-\frac{\ox^{i+1}_m-\ox^i_m}{h_m}-g(\ox^i_m,\ou^i_m)\in N(\ox^i_m;C)\;\textrm{ whenever }\;i=0,\ldots,2^m-1.
\end{equation*}
The normal cone representation mentioned in \eqref{F} ensures the existence of vectors $\eta^i_m\in\R^s_+$ for $i=0,\ldots,2^m-1$ such that conditions \eqref{ns87} and \eqref{nseta} hold for all such indices, where $\eta^i_m\in\R^s_+$ are uniquely defined due to the imposed PLICQ standing assumption. Denoting $\eta^{2^m}_m:=\xi_m$ with $\xi_m$ taken from Theorem~\ref{Thm5.2*}, we we get $\eta^i_m\in\R^s_+$ for all $i=0,\ldots,2^m$ and thus verify \eqref{nseta} for $i=2^m$. This also allows us to deduce the conditions in \eqref{nsntc} and \eqref{nsnmutx} from those in \eqref{e:5.8*} and \eqref{mutx}, respectively.\vspace*{-0.05in}

Employing now the second-order inclusion \eqref{c57ns} from Theorem~\ref{Thm6.1ns} with $x:=\ox^i_m,\;u:=\ou^i_m,\;\o:=-\frac{\ox^{i+1}_m-\ox^i_m}{h_m}$ and $w:=-\frac{1}{h_m}\lm_m\th^{iy}_m+p^{i+1}_m$ for $i=0,\ldots,2^m-1$ gives us $\gg^i_m\in\R^s$ such that
\begin{eqnarray*}
\(\frac{p^{i+1}_m-p^{i}_m}{h_m},-\frac{1}{h_m}\lm_m\th^{iu}_m-\frac{1}{h_m}\psi^i_m\)
\end{eqnarray*}
\begin{eqnarray*}
\in\partial\Big\la-\frac{1}{h_m}\lm_m\th^{iy}_m+p^{i+1}_m,g\Big\ra(\ox^i_m,\ou^i_m)+\Bigg(\begin{matrix}
\underset{j\in I_0\(p^{i+1}_m-\frac{1}{h_m}\lm_m\th^{iy}_m\)\cup I_>\(p^{i+1}_m-\frac{1}{h_m}\lm_m\th^{iy}_m\)}{\disp\sum}\gg^{ij}_m x^j_*,\\
0\end{matrix}\Bigg),
\end{eqnarray*}
with the validity of all the conditions in \eqref{nsconx}, \eqref{ns94}, and \eqref{ns93}. Implication \eqref{nseta} for $i=2^m$ follows directly from \eqref{xi} and the definition of
$\eta^{2^m}_m$. If in addition the vectors $\{x^j_*\;|\;j\in I(\ox^i_m)\}$ are linearly independent, then it follows from \eqref{ns87} and \eqref{nsconx} due to the coderivative domain formula \eqref{cod-dom} of Theorem~\ref{Thm6.1ns} that condition \eqref{ns96} is satisfied.

It remains to justify the fulfillment of the nontriviality condition \eqref{nsentc} under the imposed standing assumptions. First observe directly from \eqref{e:5.8*} that
\begin{equation}\label{nsntc}
\lm_m+\n\eta^{2^m}_m\en+\sum_{i=0}^{2^m-1}\n p^{i}_m\en+\n\psi_m\en+\|\gg_m\|\ne 0.
\end{equation}
Supposing by contraposition that \eqref{nsentc} fails gives us $\lm_m=0$, $p^{2^m}_m=0$, $\psi_m=0$, and $\gg_m=0$. Then it follows from \eqref{nsconx} that $p^i_m=0$ for all $i=0,\ldots,2^m-1$, and hence $\sum_{j=1}^s\eta^{2^mj}_m x^j_*=0$ by \eqref{nsnmutx}. This tells us that $\eta^{2^m}_m=0$ due to the assumed PLICQ. It contradicts \eqref{nsntc} and so ends the proof. $\h$\vspace*{-0.2in}

\section{Optimality Conditions for Nonsmooth Sweeping Processes}\label{optim}
\setcounter{equation}{0}\vspace*{-0.1in}

This section is the culmination of the paper. It establishes a full set of necessary optimality conditions for the original nonsmooth sweeping control problem $(P)$ formulated in \eqref{cost1}--\eqref{state} by employing the method of discrete approximations with the usage of the results obtained in the previous sections. We proceed by passing to the limit from the discrete optimality conditions from Theorem~\ref{Nonsmooth} under the standing assumptions imposed along a given relaxed $W^{1,2}\times L^2$-local minimizer of $(P)$. The reader can see that the limiting procedure is nontrivial and, besides using the assertions of Theorems~\ref{ThmStrong} and \ref{Nonsmooth}, it strongly exploits robust properties and calculations of the basic first-order and second-order generalized differential constructions discussed in Sections~\ref{tools} and \ref{2nd-order}. The obtained necessary optimality conditions extend, with different formulations and proofs in some significant points, those from \cite[Theorem~7.1]{cmn}, where it is assumed in addition that the perturbation mapping $g(x,u)$ is ${\cal C}^1$-smooth with respect to both variables and its partial Jacobian $\nabla_u g(\ox(t),\ou(t))$ has full rank on $[0,T]$ along the local minimizer in question.\vspace*{-0.03in}

Here is the main result of the paper formulated entirely in terms of the given data of $(P)$.\vspace*{-0.1in}

\begin{theorem}{\bf(necessary conditions for relaxed local minimizers).}\label{nsThm6.1*} Let $(\ox(\cdot),\ou(\cdot))$ be a relaxed $W^{1,2}\times L^2$-local minimizer for problem $(P)$ under the standing assumptions listed above. Then there exist a multiplier $\lm\ge 0$, adjoint arcs $p(\cdot)\in W^{1,2}([0,T];\R^n)$ and $q(\cdot)\in BV([0,T];\R^n)$, and a vector measure $\gg=(\gg^1,\ldots,\gg^s)\in C^*([0,T];\R^s)$ such that the following conditions are satisfied:\\[1ex]
$\bullet$ {\sc Primal arc representation}
\begin{equation}\label{ns37}
-\dot{\ox}(t)=\sum_{j=1}^s\eta^j(t)x^j_*+g\big(\ox(t),\ou(t)\big)\;\textrm{ for a.e. }\;t\in[0,T],
\end{equation}
where the functions $\eta^j(\cdot)\in L^2([0,T];\R_+)$ as $j=1,\ldots,s$ are uniquely determined by representation \eqref{ns37} for a.e.\ $t\in[0,T)$ while being well-defined at $t=T$. In fact, \eqref{ns37} holds at all $t\in[0,T]$ with $\dot{\ox}(t)$ denoting the right derivative on $[0,T)$ and postulating $\dot{\ox}(T)$ as the right-hand side of \eqref{ns37}.\\[1ex]
$\bullet$ {\sc Extended Euler-Lagrange inclusions:}
\begin{equation}\label{nsc:6.6}
\big(\dot{p}(t),-\psi(t)\big)\in\co\partial\big\la q(t),g\big\ra\big(\ox(t),\ou(t)\big)\;\mbox{ with}
\end{equation}
\begin{equation}\label{loc_max}
\psi(t)\in\co N\big(\ou(t);U\big)\;\textrm{ for a.e. }\;t\in[0,T],
\end{equation}
where $\psi(\cdot)\in L^2([0,T];\R^d)$, and where $q\colon[0,T]\to\R^{n} $ is a function of bounded variation on $[0,T]$ with its left continuous representative given, for all $t\in[0,T]$ except at most a countable subset, by
\begin{equation}\label{nsc:6.9}
q(t):=p(t)-\int_{(t,T]}\sum_{j=1}^sd\gg^j(\tau)x^j_*.
\end{equation}
$\bullet$ {\sc Maximization conditions:}
If the normal cone \eqref{nor_con} is tangentially generated
\begin{equation}\label{dua}
N(\ou(t);U)=T^*(\ou(t);U):=\big\{v\in\R^n\big|\;\la v,u\ra\le 0\;\mbox{ for all }\;u\in T\big(\ou(t);U\big)\big\}
\end{equation}
by a tangent set $T(\ou(t);U)$ associated with $U$ at $\ou(t)$, then we have the {\sc local maximization condition}
\begin{equation}\label{tan-max}
\big\la\psi(t),\ou(t)\big\ra=\max_{u\in T(\ou(t);U)}\big\la\psi(t),u\big\ra\;\textrm{ for a.e. }\;t\in[0,T]
\end{equation}
The {\sc global maximization condition}
\begin{equation}\label{max}
\big\la\psi(t),\ou(t)\big\ra=\max_{u\in U}\big\la\psi(t),u\big\ra\;\textrm{ for a.e. }\;t\in[0,T]
\end{equation}
holds if the control constraint set $U$ is convex.\\[1ex]
$\bullet$ {\sc Transversality conditions:} There exist numbers $\eta^j(T)\ge 0$ for $j\in I(\ox(T))$ such that
\begin{equation}\label{ns42}
-p(T)-\sum_{j\in I(\ox(T))}\eta^j(T)x^j_*\in\lm\partial\vph\big(\ox(T)\big)\;\mbox{ and }\;
\sum_{j\in I(\ox(T))}\eta^j(T)x^j_*\in N\big(\ox(T);C\big).
\end{equation}
$\bullet$ {\sc Complementary slackness conditions:}  We have the implications
\begin{equation}\label{ns41}
\la x^j_*,\ox(t)\ra<c_j\Longrightarrow\eta^j(t)=0\;\mbox{ and }\;\eta^j(t)>0\Longrightarrow\la x^j_*,q(t)\ra=0\;\mbox{ as }\;j=1,\ldots,s
\end{equation}
for a.e.\ $t\in[0,T)$ provided that the vectors $\{x^j_*\;|\;j\in I(\ox(t))\}$ are linearly independent for such $t$. Taking $\eta^j(T)$ from \eqref{ns41} in the agreement with \eqref{ns37}, we keep the first implication in \eqref{ns41} for $t=T$ without imposing the linear independence of $\{x^j_*\;|\;j\in I(\ox(T))\}$.\\[1ex]
$\bullet$ {\sc Measure nonatomicity conditions:} If $t\in[0,T)$ and $\la x^j_*,\ox(t)\ra<c_j$ for all $j=1,\ldots,s$, then there exists a neighborhood $V_t$ of $t$ in $[0,T)$ such that $\gg(V)=0$ for all Borel subsets $V$ of $V_t$.\\[1ex]
$\bullet$ {\sc General nontriviality conditions}: We always have
\begin{equation}\label{nse:83}
\(\lm,p,\|\gg\|\)\ne 0
\end{equation}
with the norm of $\gg$ in the space $C^*([0,T];\R^s)$. Furthermore, \eqref{nse:83}  is equivalent to $(\lm,p,q)\ne 0$ if the vectors $\{x^j_*\;|\;j\in I(\ox(t))\}$ are linearly independent on $[0,T]$.\\[1ex]
$\bullet$ {\sc Enhanced nontriviality condition}: We have $(\lm,p)\ne 0$ provided that the inequalities $\la x^j_*,\ox(t)\ra<c_j$  hold for all $t\in[0,T)$ and all indices $j=1,\ldots,s$.
\end{theorem}\vspace*{-0.07in}
{\bf Proof.} Let us split the proof of the theorem into seven steps as follows.\\[0.2ex]
{\bf Step~1:} {\em Justification of the primal arc representation}. Denote $\th^i_m:=(\th^{iy}_m,\th^{iu}_m)$ for the discrete-time functions taken from \eqref{theta} and define $\theta_m\colon[0,T)\to\R^n\times\R^d$ by
\begin{equation*}
\th_m(t):=\frac{\th^{i}_m}{h_m}\;\textrm{ for }\;t\in[t^i_m,t^{i+1}_m)\;\textrm{ and }\;i=0,\ldots,2^m-1
\end{equation*}
whenever $m\in\N$. Then we immediately arrive at
\begin{eqnarray*}
\int_0^T\n\th^y_m(t)\en^2dt&=&\sum_{i=0}^{2^m-1}\frac{\Big\|\th^{iy}_m\Big\|^2}{h_m}\le\frac{1}{h_m}\sum_{i=0}^{2^m-1}
\Big(\int_{t^i_m}^{t^{i+1}_m}\n\dot{\ox}(t)-\frac{\ox^{i+1}_m-\ox^i_m}{h_m}\en dt\Big)^2\\
&\le&\sum_{i=0}^{2^m-1}\int_{t^i_m}^{t^{i+1}_m}\Big\|\dot{\ox}(t)-\frac{\ox^{i+1}_m-\ox^i_m}{h_m}\Big\|^2dt=
\int_0^T\n\dot{\ox}(t)-\dot{\ox}_m(t)\en^2dt.
\end{eqnarray*}
Taking into account the strong convergence $(\ox_m(\cdot),\ou_m(\cdot))\to(\ox(\cdot),\ou(\cdot))$ from Theorem~\ref{ThmStrong} yields
\begin{equation}\label{c:6.14}
\int_0^T\n\th^y_m(t)\en^2dt\le\int_0^T\n\dot{\ox}(t)-\dot{\ox}_m(t)\en^2dt\to 0\;\textrm{ as }\;m\to\infty.
\end{equation}
Passing to a subsequence if needed tells us that $\th^y_m(t)\to 0$ for a.e.\ $t\in[0,T]$. Furthermore, we have
\begin{eqnarray*}
\int_0^T\Big\|\th^u_m(t)\Big\|^2dt&=&\sum_{i=0}^{2^m-1}\frac{\Big\|\th^{iu}_m\Big\|^2}{h_m}\le\frac{1}{h_m}
\sum_{i=0}^{2^m-1}\Big(\int_{t^i_m}^{t^{i+1}_m}\n\ou^i_m-\ou(t)\en dt\Big)^2\\
&\le&\sum_{i=0}^{2^m-1}\int_{t^i_m}^{t^{i+1}_m}\n\ou^i_m-\ou(t)\en^2dt=\int_0^T\n\ou_m(t)-\ou(t)\en^2dt,
\end{eqnarray*}
and thus deduce from Theorem~\ref{ThmStrong} that
\begin{equation}\label{c:6.14b}
\int_0^T\n\th^u_m(t)\en^2dt\le\int_0^T\n\ou_m(t)-\ou(t)\en^2dt\to 0,\;\mbox{ and so }\;\th^u_m(t)\to 0\;\mbox{ a.e. on }\;[0,T]
\end{equation}
without relabeling. To proceed further, observe that the assumed PLICQ condition is robust with respect to perturbations of the initial point, and hence it also holds along the approximating vectors $\ox^i_m$ as $m\to\infty$. This allows us to construct the piecewise constant functions $\eta_m\colon[0,T)\to\R^s_+$ by
\begin{equation*}
\eta_m (t):=\eta^{i}_m\;\mbox{ for }\;t\in[t^i_m,t^{i+1}_m),\quad i=0,\ldots,2^m-1,
\end{equation*}
where $\eta^{i}_m\in\R^s_+$ are uniquely determined by Theorem~\ref{Nonsmooth}. It follows from \eqref{ns87} that
\begin{equation}\label{nsc:51}
-\dot{\ox}_m(t)=\sum_{j=1}^s\eta^j_m(t)x^j_*+g\big(\ox_m(t^i_m),\ou_m(t^i_m)\big)-r_m(t^i_m)\rho_m(t^i_m)\;\textrm{ for all }\;t\in(t^i_m,t^{i+1}_m).
\end{equation}
Passing to the limit in \eqref{nsc:51} as $m\to\infty$ and using Theorem~\ref{ThmStrong} yield
\begin{equation*}
-\dot{\ox}(t)\in N\big(\ox(t);C\big)+g\big(\ox(t),\ou(t)\big)\;\mbox{ for a.e. }\;t\in[0,T].
\end{equation*}
Remembering the representation of the normal cone in \eqref{F} and taking into account the imposed PLICQ property along $\ox(t)$ allow us to conclude that the mapping
\begin{equation}\label{eta-mes}
(\lm^j)_{j\in I(\ox(t))}\colon[0,T]\mapsto\dot{\ox}(t)+g\big(\ox(t),\ou(t)\big)+\Big\{\sum_{j\in I(\ox(t))}\lm^j x^j_*\;\Big|\;\lm^j\ge 0\Big\}
\end{equation}
is single-valued on $[0,T)$. Denote by $\lm^j(t)$ for $j\in I(\ox(t))$ the functions from \eqref{eta-mes} and then define $\eta^j\colon[0,T]\to\R_+$ for all $j\in\{1,\ldots,s\}$ by
\begin{equation*}
\eta^j(t):=\left\{\begin{array}{ll}
\lm^j(t)&\mbox{for }\;j\in I(\ox(t)),\\
0&\mbox{otherwise.}
\end{array}\right.
\end{equation*}
and deduce from \eqref{eta-mes}, due to $\dot{\ox}(\cdot)\in L^2([0,T];\R^n)$ and (H3), that each $\eta^j(t)$ as $j\in\{1,\ldots,s\}$ belongs the space $L^2([0,T];\R_+)$. This clearly verifies the claimed primal arc representation \eqref{ns37} for a.e.\ $t\in[0,T)$ with the remark on the right derivative of $\dot{\ox}(t)$ therein. To define $\eta^j(t)$ at the endpoint $t=T$, take $\eta_m(T):=\eta^{2^m}_m$ from the optimality conditions for discrete approximations in Theorem~\ref{Nonsmooth} and deduce from the nontriviality conditions in \eqref{nsntc} after their normalization that the sequence $\{\eta^{2^m}_m\}$ converges, along a subsequence, to some vector $(\eta^{1}(T),\ldots,\eta^{2^m}(T))$, which is used in what follows.\\[1ex]
{\bf Step~2:} {\em Construction of approximating arcs}. Let us first consider the piecewise linear extensions of the discrete-time functions $p^i_m$ from Theorem~\ref{Nonsmooth} to the continuous-time interval $[0,T]$ and define $q_m(t)$ by $q_m(t^i_m):=p^i_m$ for $i=0,\ldots,2^m$. Then $\gg_m(t)$ and $\psi_m(t)$ on $[0,T]$ are constructed by
\begin{equation}\label{const1}
\gg_m(t):=\gg^i_m\;\mbox{ and }\;\psi_m(t):=\frac{1}{h_m}\psi^i_m\;\textrm{ for }\;t\in[t^i_m,t^{i+1}_m)\;\textrm{ and }\;i=0,\ldots,2^m-1
\end{equation}
with $\gg_m(T):=0$ and $\psi_m(T):=0$. Define further
\begin{equation*}
\nu_m(t):=\max\big\{t^i_m\;\big|\;t^i_m\le t,\;0\le i\le 2^m-1\big\}\;\textrm{ for all }\;t\in[0,T]\;\mbox{ and }\;m\in\N.
\end{equation*}
It follows from the inclusions in \eqref{nsconx} that
\begin{eqnarray}\label{nsconx'}
\begin{array}{ll}
\disp\big(\dot{q}_m(t),-\lm_m\th^{u}_m(t)-\psi_m(t)\big)&\in\partial\Big\la-\lm_m\th^{y}_m(t)+q_m\big(\nu_m(t)+h_m\big),g\Big\ra\big(\ox_m(\nu_m(t)),\ou_m(\nu_m(t))\big)\\\\
&+\disp\(\sum_{j\in I_0\(p^{i+1}_m-\frac{1}{h_m}\lm_m\th^{iy}_m\)\cup I_>\(p^{i+1}_m-\frac{1}{h_m}\lm_m\th^{iy}_m\)}\gg^{j}_m(t)x^j_*,\;0\)
\end{array}
\end{eqnarray}
for every $t\in (t^i_m,t^{i+1}_m),\;i=0,\ldots,2^m-1$. Define now $p_m(\cdot)$ on $[0,T]$ by setting
\begin{equation}\label{nsc:6.29}
p_m(t):=q_m(t)+\int_t^T\Big(\sum_{j=1}^s\gg^j_m(\tau)x^j_\ast\Big)d\tau\;\textrm{ for every }\;t\in[0,T],
\end{equation}
which gives us $p_m(T)=q_m(T)$ and the differential relations
\begin{equation}\label{nsc:6.30}
\dot{p}_m(t)=\dot{q}_m(t)-\sum_{j=1}^s\gg^j_m(t)x^j_\ast\;\textrm{ for a.e. }\;t\in[0,T].
\end{equation}
It follows from $\[j\notin I_0\(p^{i+1}_m-\frac{1}{h_m}\lm_m\th^{iy}_m\)\cup I_>\(p^{i+1}_m-\frac{1}{h_m}\lm_m\th^{iy}_m\)\Longrightarrow\gg^{ij}_m=0\]$ in $\eqref{ns93}$ that
\begin{equation}\label{sum}
\sum_{j=1}^s\gg^j_m(\tau)x^j_\ast=\sum_{j\in I_0\(p^{i+1}_m-\frac{1}{h_m}\lm_m\th^{iy}_m\)\cup I_>\(p^{i+1}_m-\frac{1}{h_m}\lm_m\th^{iy}_m\)}\gg^{j}_m(t)x^j_\ast.
\end{equation}
Using this together with \eqref{nsconx'}, \eqref{nsc:6.30}, and \eqref{sum} yields
\begin{equation}\label{nsc:59}
\(\dot{p}_m(t),-\lm_m\th^{u}_m(t)-\psi_m(t)\)\in\partial\la-\lm_m\th^{y}_m(t)+q_m(\nu_m(t)+h_m),g\ra\big(\ox_m(\nu_m(t)),\ou_m(\nu_m(t))\big)
\end{equation}
for every $t\in(t^i_m,t^{i+1}_m),\;i=0,\ldots,m-1$. Next we construct the vector measures $\gg^{mes}_m$ on $[0,T]$ by
\begin{equation}\label{nsc:6.34}
\underset{B}{\int}d\gg^{mes}_m:=\underset{B}{\int}\sum_{i=0}^{2^m-1}\dfrac{1}{h_m}\gg_m(t)\mathbbm{1}_{I^i_m}(t)dt,
\end{equation}
for every Borel subset $B\subset[0,T]$, and where $\mathbbm{1}_{I^i_m}$ stands for the characteristic function of the set $I^i_m$. Dropping from now on the index ``$mes$" in the measure notation and observing that all the expressions in the statement of Theorem~\ref{Thm5.2*} are positively homogeneous of degree 1 with respect to $(\lm_m,\,p_m,\,\gg_m,\,\psi_m)$, we normalize the nontriviality condition \eqref{nsentc} with the usage of the constructions in \eqref{const1} by
\begin{equation}\label{nsc:6.35}
\lm_m+\n p_m(T)\en+\sum_{j=1}^s\sum_{i=0}^{2^m-1}\big|\gg^{ij}_m\big|+\int_0^T\n\psi_m(t)\en dt=1,\quad m\in\N.
\end{equation}
This tells us that all the sequences in \eqref{nsc:6.35} are uniformly bounded, which is employed in what follows.\\[1ex]
{\bf Step~3:} {\em Verification of the extended Euler-Lagrange and maximization conditions.} Passing to a subsequence if needed, we get from \eqref{nsc:6.35} that $\lm_m\to\lm$ as $m\to\infty$ for some $\lm\ge 0$. Let us now verify that the sequence $\{ p^{0}_m,\ldots,p^{2^m}_m\}_{m\in\N}$ is bounded. Indeed, it follows from \eqref{nsconx} that
\begin{equation*}\label{nsc-p}
\(\frac{p^{i+1}_m-p^{i}_m}{h_m}-\sum_{j =1 }^{s}\gg^{ij}_mx^j_*,\;-\frac{1}{h_m}\lm_m\th^{iu}_m-\frac{1}{h_m}\psi^i_m\)\in\partial\Big\la-\frac{1}{h_m}\lm_m\th^{iy}_m+p^{i+1}_m,g\Big\ra (\ox^i_m,\ou^i_m)
\end{equation*}
for all $i=0,\ldots,2^m-1$. The latter implies, by taking into account that the subgradient sets for locally Lipschitzian functions are bounded (see, e.g., \cite[Theorem~1.22]{m18}) and employing \eqref{c:6.14} together with \eqref{nsc:6.35}, that the sets $\partial\big\la-\frac{1}{h_m}\lm_m\th^{iy}_m+p^{i+1}_m,g\big\ra (\ox^i_m,\ou^i_m)$ are uniformly bounded for all $i=0,\ldots,2^m-1$ and $m\in\N$. Thus there exists a constant $M_1>0$ ensuring that
\begin{equation}\label{M1}
\Big\|\Big(\frac{p^{i+1}_m-p^{i}_m}{h_m}-\sum_{j=1}^{s}\gg^{ij}_mx^j_*,\;-\frac{1}{h_m}\lm_m\th^{iu}_m-\frac{1}{h_m}\psi^i_m\Big)\Big\|\le M_1
\end{equation}
for all these indices. Using the sum norm in $\R^{n+d}$, we get the estimate
\begin{equation}\label{pq}
\Big\|\frac{p^{i+1}_m-p^{i}_m}{h_m}-\sum_{j=1}^{s}\gg^{ij}_mx^j_*\Big\|\le M_1,
\end{equation}
which implies in turn the relationships
\begin{equation*}
\big\|p^{i}_m\big\|\le\big\|p^{i+1}_m\big\|+\n h_m\en M_1+\Big\|h_m\sum_{j=1}^{s}\gg^{ij}_mx^j_*\Big\|=\n p^{i+1}_m\en+h_mM_1+h_m\Big\|\sum_{j=1}^{s}\gg^{ij}_mx^j_*\Big\|
\end{equation*}
for all $i=0,\ldots,2^m-1$. For each $m\in\N$ consider the quantity
\begin{equation*}
A^i_m:=h_mM_1+h_m\Big\|\sum_{j=1}^{s}\gg^{ij}_mx^j_*\Big\|,\quad i=0,\ldots,2^m,
\end{equation*}
where $h_m M_1\to 0$ as $m\to\infty$. On the other hand, it follows from \eqref{nsc:6.35} that
\begin{equation}\label{2ndterm}
\sum_{i=0}^{2^m-1}h_m\Big\|\sum_{j=1}^s\gg^{ij}_mx^j_*\Big\|=\int_0^T\Big\|\sum_{j=1}^s\gg^j_m(t)x^j_*\Big\|dt\le 1.
\end{equation}
This yields $\sum_{i=0}^{2^m-1}A^i_m\le M_2$ for some $M_2>0$. Combining it with the estimate for $\|p^i_m\|$ ensures that
\begin{equation*}\label{t:67}
\n p^{i}_m\en\le\n p^{i+1}_m\en+A^i_m\;\mbox{ for }\;i=0,\ldots,2^m-1.
\end{equation*}
Using this step by step, we arrive at
\begin{eqnarray*}
\n p^{i}_m\en\le\n p^{2^m}_m\en+\sum_{j=i}^{2^m-1}A^j_m&\le 1+\disp\sum_{i=0}^{2^m-1}A^i_m\le 1+M_2\;\mbox{ as }\;i=0,\ldots,2^m-1,
\end{eqnarray*}
which verifies the boundedness of $\{p^{i}_m\}_{0\le i\le 2^m}$. To deal now with the functions $q_m(\cdot)$, derive from their construction in Step~2 and the estimate in \eqref{pq} that
\begin{equation*}
\bigg\|\sum_{i=0}^{2^m-1}\big\|q_m(t^{i+1}_m)-q_m(t^i_m)\big\|-\int_0^T\Big\|\sum_{j=1}^s\gg^{j}_m(t)x^j_*\Big\|dt\bigg\|\le h_m M_1,
\end{equation*}
which tells us therefore that
\begin{equation}\label{t:68}
\sum_{i=0}^{2^m-1}\n q_m(t^{i+1}_m)-q_m(t^i_m)\en\le h_m M_1+\int_0^T\Big\|\sum_{j=1}^s\gg^{j}_m(t)x^j_*\Big\|dt.
\end{equation}
It follows from \eqref{2ndterm}, \eqref{t:68}, and the construction of $q_m(t)$ on $[0,T]$ that these functions are of uniformly bounded variations on this interval. Due to the obvious inequalities
\begin{equation*}
2\n q_m(t)\en-\n q_m(0)\en-\n q_m(T)\en\le\n q_m(t)-q_m(0)\en+\n q_m(T)-q_m(t)\en\le\textrm{var}(q_m;[0,T]),\quad t\in[0,T],
\end{equation*}
and the boundedness of $\{q_m(0)\}$ and $\{q_m(T)\}$, we get the uniform boundedness of the functions $q_m(\cdot)$ on $[0,T]$. Then Helly's selection theorem allows us to find a function of bounded variation $q(\cdot)$ on $[0,T]$ such that $q_m(t)\to q(t)$ as $m\to\infty$ for each $t\in[0,T]$. Employing \eqref{nsc:6.35} and the measure construction in \eqref{nsc:6.34} tell us that the sequence $\{\gg_m\}$ is bounded in $C^*([0,T];\R^s)$. It follows from the weak* sequential compactness of the unit balls in this space that there are measures $\gg\in C^*([0,T];\R^s)$ such that $\{\gg_m\}$ weak* converges to $\gg$ in $C^*([0,T];\R^s)$ along some subsequence; see \cite[Proposition~3.21 and Theorem~3.23]{AFP}.\vspace*{-0.05in}

Furthermore, we have the convergence
\begin{equation*}
\Big\|\int_t^T\sum_{j=1}^s\gg^j_m(\tau)x^j_\ast d\tau-\int_{(t,T]}\sum_{j=1}^sd\gg^j(\tau)x^j_*\Big\|\to 0\;\textrm{ as }\;m\to\infty
\end{equation*}
for all $t\in [0,T]$ except a countable subset of $[0,T]$ by the weak* convergence of the measures $\gg_m$ to $\gg$ in the space $C^*([0,T];\R^n)$; see \cite[p.\ 325]{v} for similar arguments. Hence
\begin{equation*}
\int^T_t\sum_{j=1}^s\gg^j_m(\tau)x^j_\ast d\tau\to\int_{(t,T]}\sum_{j=1}^sd\gg^j(\tau)x^j_*\;\textrm{ on }\;[0,T]\;\textrm{ as }\;m\to\infty,
\end{equation*}
and thus we arrive at \eqref{nsc:6.9} by passing to the limit in \eqref{nsc:6.30} as $m\to\infty$.

Turning now to $\psi_m(\cdot)$, deduce from \eqref{const1} and \eqref{M1} that
\begin{equation*}
\big|\,\|\psi_m(t)\|-\lm_m\|\th^u_m(t)\|\,\big|\le\n-\psi_m(t)-\lm_m\th^u_m(t)\en\le M_1,
\end{equation*}
which readily brings us to the estimate
\begin{equation}\label{psi-bound}
\n\psi_m(t)\en\le\lm_m\n\th^u_m(t)\en+M_1\;\textrm{ for all }\;t\in[t^i_m,t^{i+1}_m).
\end{equation}
Taking into account \eqref{nsc:6.35} and using the Cauchy-Schwarz inequality yield
\begin{equation*}
\begin{aligned}
\int^{t^{i+1}_m}_{t^i_m}\n\psi_m(t)\en^2 dt\le 2\[\int^{t^{i+1}_m}_{t^i_m}\n\th^u_m(t)\en^2dt+M_1^2h_m\],
\end{aligned}
\end{equation*}
and therefore we get the boundedness of $\{\psi_m(\cdot)\}$ in $L^2([0,T];\R^d)$:
\begin{equation*}
\begin{aligned}
\int^T_0\n\psi_m(t)\en^2dt\le 2\[\underbrace{\int^T_0\n\th^u_m(t)\en^2dt}_{\dn 0\;\;\mbox{as}\;\;m\to\infty}+2^m\dfrac{T}{2^m}M_1^2\].
\end{aligned}
\end{equation*}
Thus there exists a subsequence of $\{\psi_m(\cdot)\}$, which weakly converges to some function $\psi(\cdot)\in L^2([0,T];\R^d)$ in this space. Combining further the obtained uniform boundedness of $q_m(\cdot)$ on $[0,T]$ with \eqref{nsc:59} and \eqref{nsc:6.35} allows us to conclude that the sequence of the derivatives $\{\dot{p}_m(\cdot)\}$ is bounded in $L^2([0,T];\R^n)$ and thus weakly compact in this space. It tells us that a subsequence of $\{\dot{p}_m(\cdot)\}$ weakly converges to some function $v(\cdot)\in L^2([0,T];\R^n)$ as $m\to\infty$, without relabeling. Denoting $p(T):=q(T)$, where $q(T)$ is constructed above, we define the adjoint arc
\begin{equation*}
p(t):=p(T)+\int_{T}^t v(\tau)d\tau\;\mbox{ for all }\;t\in[0,T).
\end{equation*}
This ensures that $\dot p(\cdot)=v(t)$ for a.e.\ $t\in[0,T]$ and $\{p_m(\cdot)\}$ converges to $p(\cdot)$ weakly in $W^{1,2}([0,T];\R^n)$.\vspace*{-0.03in}

Applying the classical Mazur theorem to the sequence of pairs $\{(\dot{p}_m(\cdot),\psi_m(\cdot))\}$ gives us a sequence of convex combinations of $\{(\dot{p}_m(\cdot),\psi_m(\cdot))\}$, which converges to $(\dot p(\cdot),\psi(\cdot))$ strongly in $L^2([0,T];\R^n)\times L^2([0,T];\R^d)$ and thus pointwise for a.e.\ $t\in[0,T]$ along a subsequence. Passing now to the pointwise limit in \eqref{nsc:59} along this subsequence of convex combinations with taking into account that $\(\th^u_m(t),\th^y_m(t)\)\to(0,0)$ for a.e.\ $t\in[0,T]$ as proved in \eqref{c:6.14} and \eqref{c:6.14b}, we arrive at the Euler-Lagrange inclusion \eqref{nsc:6.6} with the function $\psi(\cdot)$ determined above. To verify finally that this function satisfies \eqref{loc_max}, recall from Theorem~\ref{Nonsmooth} that $\psi^{i}_m\in N(\ou^i_m;U)$ for $i=0,\ldots,2^m-1$. Remembering the construction of $\psi_m(\cdot)$ in \eqref{const1}, the piecewise constant extension of $\ou^i_m$ to $[0,T]$, and the conic structure of $N(\cdot;U)$ imply that
\begin{equation}\label{coN}
\psi_m(t)\in N\big(\ou_m(t);U\big)\;\textrm{ for all }\;t\in[t^i_m,t^{i+1}_m)\;\textrm{ and }\;i=0,\ldots,2^m-1.
\end{equation}
Then we can pass to the pointwise limit in \eqref{coN} along a subsequence of $m\to\infty$ by employing the strong $L^2$-convergence of $\ou_m(\cdot)\to\ou(\cdot)$ from Theorem~\ref{ThmStrong}, the robustness of the normal cone \eqref{nor_con} with respect to perturbations of the initial point, the strong $L^2$-convergence of convex combinations of $\psi_m(\cdot)$ to $\psi(\cdot)$, and the boundedness of $\psi(\cdot)$ on $[0,T]$ due to \eqref{nsc:6.6} under (H3) and $q(\cdot)\in BV([0,T];\R^n)$. In this way we justify the claimed inclusion \eqref{loc_max} by using the normal cone convexification.\vspace*{-0.03in}

To finish the proof in this step, it remains to verify the local and global maximization conditions in \eqref{tan-max} and \eqref{max}, respectively. Note that the duality correspondence in \eqref{dua} generated by any tangent set $T(\ou(t);U)$ always yields the convexity of $N(\ou(t);\O)$. Thus the local maximization condition \eqref{tan-max} follows directly from \eqref{loc_max} and \eqref{dua}. If $U$ is convex, then the normal cone \eqref{nor_con} and hence its convexification in \eqref{loc_max} reduce to the normal cone of convex analysis \eqref{NC}. Thus the global maximization condition \eqref{max} is an immediate consequence of \eqref{loc_max} and the structure of \eqref{NC}.\\[1ex]
{\bf Step~4:} {\em Verification of the transversality conditions.} Recalling the definition of $(\eta^{1}(T),\ldots,\eta^{2^m}(T))$ as a limiting point of $\eta_m(T):=\eta^{2^m}_m$ as $m\to\infty$ in Step~1 of the proof, we get from the discrete transversality conditions \eqref{nsnmutx} and the normal cone representation in \eqref{F} that
\begin{equation}\label{c:72}
-p^{2^m}_m-\lm_m\vth^{2^m}_m=\sum_{j=1}^s\eta^{2^mj}_m x^{2^m}_\ast=\sum_{j\in I(\ox^{2^m}_m)}\eta^{2^mj}_m x^{2^m}_\ast\in N(\ox^{2^m}_m;C),
\end{equation}
where $\eta^{2^mj}_m:=0$ for $j\in\{1,\ldots,s\}\setminus I(\ox^{2^m}_m)$. Denote $\zeta_m:=\sum_{j\in I(\ox^{2^m}_m)}\eta^{2^mj}_m x^{2^m}_\ast$ and deduce from \eqref{nsnmutx} and \eqref{c:72} due to the boundedness of $\lm_m$ by \eqref{nsc:6.35}, the convergence of $\{\ox^{2^m}_m\}$ and $\{p^{2^m}_m\}$, and the boundedness of the subdifferential of the locally Lipschitzian function $\ph$ that a subsequence of $\{\zeta_m\}$ converges to some $\zeta\in\R^n$. Using then the robustness of the normal cone in \eqref{c:72}, the convergence of $\ox^{2^m}_m\to\ox(T)$, and the inclusion $I(\ox^{2^m}_m)\subset I(\ox(T))$ for all large $m$, we get $\zeta\in N(\ox(T);C)$. It follows from \eqref{nsnmutx} that
\begin{equation}\label{trans-dis}
-p^{2^m}_m-\zeta_m\in\lm_m\partial\vph(\ox^{2^m}_m)\;\textrm{ for all }\;m\in\N.
\end{equation}
This allows us derive both transversality inclusions in \eqref{ns42} by the passage to the limit in \eqref{trans-dis} as $m\to\infty$ with taking into account the robustness of the subdifferential mapping therein.\\[1ex]
{\bf Step~5:} {\em Verification of the complementarity conditions}. It follows from \eqref{ns37} and \eqref{nsc:51} due to the constructions of Step~1 that we have
\begin{equation*}
\dot{\ox}(t)-\dot{\ox}_m(t)=\sum_{j=1}^s\big[\eta^j_m(t)-\eta^j(t)\big]x^j_*-g\big(\ox(t),\ou(t)\big)+g\big(\ox_m(t^i_m),\ou_m(t^i_m)\big)-r_m(t^i_m)\rho_m(t^i_m)
\end{equation*}
for $t\in(t^i_m,t^{i+1}_m)$ and $i=0,\ldots,2^m-1$. It gives us the estimate
\begin{equation*}
\Big\|\sum_{j=1}^s\big[\eta^j(t)-\eta^j_m(t)\big]x^j_*\Big\|_{L^2}\le\n\dot{\ox}_m(t)-\dot{\ox}(t)\en_{L^2}+\n g\big(\ox_m(t),\ou_m(t)\big)-g\big(\ox(t),\ou(t)\big)\en_{L^2}+r_m(t^i_m)
\end{equation*}
whenever $t\in(t^i_m,t^{i+1}_m)$. Using Theorem~\ref{ThmStrong} and passing to the limit in the above inequality, we get
\begin{equation*}
\sum_{j\in I(\ox(t))}\big[\eta^j(t)-\eta^j_m(t)\big]x^j_*\to 0\;\textrm{ as }\;m\to\infty\;\textrm{ for a.e. }\;t\in[0,T).
\end{equation*}
It clearly yields the first complementarity condition in \eqref{ns41} for a.e.\ $t\in[0,T)$ by the assumed linear independence of $\{x^j_*\;|\;j\in I(\ox(t))\}$ for such $t$. The fulfillment of this complementary slackness condition for $t=T$ follows from its discrete counterpart in \eqref{nseta} of Theorem~\ref{Nonsmooth} without imposing the additional linear independence assumption on the vectors $\{x^j_*\;|\;j\in I(\ox(T))\}$ due to the endpoint convergence $\ox^{2m}_m\to\ox(T)$ and $\eta^{2m}_m\to\eta(T)$ as $m\to\infty$ discussed above in Steps~1 and 4.\vspace*{-0.03in}

To justify the second qualification condition in \eqref{ns41} for a.e.\ $t\in[0,T)$ under the linear independence assumption on $\{x^j_*\;|\;j\in I(\ox(t))\}$ for the corresponding $t$, we employ the discrete ones in \eqref{ns96} valid under the imposed LICQ assumption due to its robustness when $\ox^i_m\to\ox(t)$. Recalling the constructions of $\eta^m(t),\;\nu_m(t),\;q_m(t)$, and $\th_m(t)$ from Steps~1 and 2, we rewrite \eqref{ns96} in the form
\begin{equation}\label{complem}
\eta^{j}_m(t)>0\Longrightarrow\big[\big\la x^j_*,q_m\big(\nu_m(t)+h_m\big)-\lm_m\th^{y}_m(t)\big\ra=0\big]\;\mbox{ for }\;j=1,\ldots,s\;\mbox{ and }\;t\in[0,T).
\end{equation}
Taking now into account that $\eta_m(t)\to\eta(t)$, $q_m(t)\to q(t)$, and $\th^y_m(t)\to 0$ for a.e.\ $t\in[0,T)$ as $m\to\infty$ along a subsequence, we obtain the claimed implication in \eqref{ns41} by passing to the limit in \eqref{complem}.\\[1ex]
{\bf Step~6:} {\em Verification of the measure nonatomicity conditions.} The proofs in this step and largely in Step~7 are similar to the corresponding arguments in \cite[Theorem~7.1]{cmn} (with some valuable modifications in Step~7), but we present them here for the reader convenience. To verify the measure nonatomicity condition, pick any $t\in[0,T)$ with $\la x^j_*,\ox(t)\ra<c_j$ for all $j=1,\ldots,s$ and by using the continuity of $\ox(\cdot)$ find a neighborhood $V_t$ of $t$ such that $\la x^j_*,\ox(\tau)\ra<c_j$ when $\tau\in V_t$ and $j=1,\ldots,s$. Then we employ Theorem~\ref{ThmStrong} to get $\la x^j_\ast,\ox_m(t^i_m)\ra<c_j$ whenever $t^i_m\in V_t$ for all $j=1,\ldots,s$ and all large $m\in\N$. It follows from \eqref{ns94} that $\gg_m(t)=0$ on any Borel subset $V$ of $V_t$. Thus
\begin{equation}\label{nonatom}
\|\gg_m\|(V)=\disp\int_Vd\|\gg_m\|=\int_V\|\gg_m(t)\|dt=0
\end{equation}
by the construction of $\gg_m(\cdot)$ in \eqref{nsc:6.34}. The passage to the limit in \eqref{nonatom} with taking into account the measure convergence obtained in Step~3 gives us $\gg(V)=0$ for the limiting measure and hence verifies the claimed measure nonatomicity condition.\\[1ex]
{\bf Step~7:} {\em Verification of the nontriviality conditions.} Let us first justify the general nontriviality condition \eqref{nse:83} without any additional assumptions. Suppose by contraposition that $\lm=0$, $p(t)=0$ for all $t\in[0,T]$, and $\|\gg\|=0$. Dealing with the right continuous representative of $q(\cdot)$, we obtain that $q(t)=0$ due to \eqref{nsc:6.9}. Thus $\lm_m\to 0$ and $p_m(t)\to 0$ for all $t\in[0,T]$ by the assumed negation of the nontriviality condition. Due to \eqref{nsc:6.29} and the convergence result taken from \cite[p.\ 325]{v}, we get that
\begin{equation*}
\underset{m\to\infty}{\lim}q_m(t)=\underset{m\to\infty}{\lim}\Big(p_m(t)-\int_t^T\sum_{j=1}^s\gg^j_m(\tau)x^j_*d\tau\Big)=\underset{m\to\infty}{\lim}p_m(t)-\int_t^T\sum_{j=1}^s\gg(\tau)x^j_*d\tau=0.
\end{equation*}
This allows us to deduce from \eqref{nsconx'} that $\psi_m(t)\to 0$ for a.e.\ $t\in[0,T]$. Using it together with estimate \eqref{psi-bound} tells us by the dominated convergence theorem that
\begin{equation*}
\int^T_0\|\psi_m(t)\|dt\to 0,\;\mbox{ and so }\;\sum_{j=1}^s\sum_{i=0}^{2^m-1}\;\big|\gg^{ij}_m\big|\to 1\;\mbox{ as }\;m\to\infty
\end{equation*}
due to \eqref{nsc:6.35}. By using the Jordan measure decompositions
\begin{equation*}
\gg_m=(\gg_m)^+-(\gg_m)^-\;\mbox{ and }\;\gg=\gg^{+}-\gg^{-}
\end{equation*}
and recalling the separability of the space $C^*([0,T];\R^s)$, we find a subsequence of measures $\{\gg_m\}$ with
\begin{equation*}
\{(\gg_m)^+\}\overset{w^*}{\to}\gg^+\;\textrm{ and }\;\{(\gg_m)^-\}\overset{w^*}{\to}\gg^-\;\textrm{ in }\;C^*\([0,T];\R^s\)\;\mbox{ as }\;m\to\infty,
\end{equation*}
where $w^*$ signifies the weak$^*$ topology. For each $m\in\N$ define the mapping $\alpha_m\colon[0,T]\to\R^s$ by
\begin{equation*}
\alpha^i_m(t):=\dfrac{\gg^i_m(t)}{|\gg^i_m(t)|}\;\textrm{ if }\;\gg^i_m(t)\ne 0\;\textrm{ and }\;\alpha^i_m(t):=0\;
\textrm{ if }\;\gg^i_m(t)=0\;\mbox{ for all }\;i=1,\ldots,s.
\end{equation*}
Since the functions $\alpha_m(\cdot)$ are clearly measurable and uniformly bounded on $[0,T]$, the application of \cite[Proposition~9.2.1]{v} (where our index $m$ corresponds to the index $i$ in that result) with $A=A_m:=[-1,1]^s$ for all $m\in\N$ ensures the existence of Borel measurable functions $\alpha^+,\,\alpha^-:[0,T]\to\R^s$ satisfying
\begin{equation*}
\{\alpha^i_m(\gg^i_m)^+\}\overset{w^*}{\to}(\alpha^+)^i(\gg^+)^i\;\textrm{ and }\;\{\alpha^i_m(\gg^i_m)^-\}\overset{w^*}{\to}(\alpha^-)^i(\gg^-)^i\;\mbox{ as }\;m\to\infty\;\mbox{ for all }\;i=1,\ldots,s.
\end{equation*}
Denoting $\alpha d\gamma:=(\alpha^i d\gamma^i,\ldots,\alpha^s d\gamma^s)$, we get the following relationships:
\begin{equation*}
\begin{aligned}
&\n\int_{[0,T]\backslash S}\alpha^+(t)d\gg^+(t)-\int_{[0,T]\backslash S}\alpha^-(t)d\gg^-(t)\en=\lim_{m\to\infty}\n\int_{[0,T]\backslash S}\alpha_m(t)d\(\gg_m\)^+(t)-\int_{[0,T]\backslash S}\alpha_m(t)d\(\gg_m\)^-(t)\en\\
&=\lim_{m\to\infty}\n\int_{[0,T]\backslash S}\alpha_m(t)d\gg_m(t)\en=\lim_{m\to\infty}\n\int_{[0,T]\backslash S}\(\alpha_m^1(t)d\gg_m^{1}(t),\ldots,\alpha_m^s(t)d\gg_m^{s}(t)\)\en\\
&=\lim_{m\to\infty}\n\(\sum^{2^m-1}_{i=0}\left|\gg_m^{i1}\right|,\ldots,\sum^{2^m-1}_{i=0}\left|\gg_m^{is}\right|\)\en=\lim_{m\to\infty}\sqrt{\sum^s_{j=1}\[\sum^{2^m-1}_{i=0}\left|\gg_m^{ij}
\right|\]^2}\ge\lim_{m\to\infty}\dfrac{1}{\sqrt s}\sum^s_{j=1}\sum^{2^m-1}_{i=0}\left|\gg_m^{ij}\right|=\dfrac{1}{\sqrt s}>0,
\end{aligned}
\end{equation*}
where $S\subset[0,T]$ is a countable set. Moreover, we also have the estimates
\begin{equation*}
\begin{aligned}
&\n\int_{[0,T]\backslash S}\alpha^+(t)d\gg^+(t)-\int_{[0,T]\backslash S}\alpha^-(t)d\gg^-(t)\en\le\n\int_{[0,T]\backslash S}\alpha^+(t)d\gg^+(t)\en+\n\int_{[0,T]\backslash S}\alpha^-(t)d\gg^-(t)\en\\
&\le\int_{[0,T]\backslash S}d\|\gg^+(t)\|+\int_{[0,T]\backslash S}d\|\gg^-(t)\|\le\n\gg^+\en+\n\gg^-\en=\n\gg\en.
\end{aligned}
\end{equation*}
Combining all the above tells us that $\n\gg\en>0$, which contradicts the contraposition assumption $\n\gg\en=0$ and thus verifies \eqref{nse:83}. If we additionally impose the linear independence of $\{x_*^j\;|\;j\in I(\ox(t))\}$ on $[0,T]$, than the obtained general nontriviality conditions is clearly equivalent to $(\lm,p,q)\ne 0$ due to \eqref{nsc:6.9}.\vspace*{-0.03in}

To complete the proof of the theorem, we need to verify the enhanced nontriviality condition $(\lm,p)\ne 0$ under the interiority assumptions imposed therein. If it is not the case, then $\lm=0$ and $p(t)\equiv 0$ on $[0,T]$ in spite of $\la x^j_*,\ox(t)\ra<c_j$ for all $t\in[0,T)$ and $j=1,\ldots,s$. Assuming the latter, we deduce for the discrete complementarity slackness condition in \eqref{ns94} and the arguments in Step~5 with the usage of \eqref{nonatom} that $\|\gg\|=0$. This contradicts \eqref{nse:83} and thus finishes the proof. $\h$\vspace*{-0.2in}

\section{Applications to Robotics}\label{applic}
\setcounter{equation}{0}\vspace*{-0.1in}

In this concluding section of the paper, we present applications of the obtained necessary optimality conditions for the nonsmooth sweeping control problem $(P)$ to some dynamical model appearing in robotics. The uncontrolled dynamics of this model, known as the {\em mobile robot model with obstacles}, was described as a sweeping process in \cite{HB}. The recent paper \cite{cmn1} contains an optimal control formalization of this model as a perturbed controlled sweeping process. Due to the scope of theoretical results for optimization of controlled sweeping processes available at that time, only a smooth version of the mobile robot model has been investigated in \cite{cmn1}. However, a more realistic version of this model requires the usage of nonsmooth controlled perturbations, which would allow the controllers to adequately react on the sudden change of the robot velocity at the contact time with an obstacle. Indeed, the velocity of the robot should be adjusted in order to keep the distance from other obstacles by using some control actions in the velocity term that is also treated as perturbations. Moreover, since the robot velocity can be changed suddenly by other obstacles, the graph of a perturbation function has some abrupt bends, cusps, and/or corners, which lead us to considering nonsmooth perturbations. Observe further that in real-life models, motions of the robots are usually accompanied by dry friction forces; see, e.g., \cite{Cher}. By Coulomb's law, dry friction forces are defined by nonsmooth functions of the velocity of the moving robot, and thus the system dynamics is naturally described by differential inclusions with nonsmooth right-hand sides. This also motivates our interest to study mobile robot models with nonsmooth perturbations.\vspace*{-0.03in}

\begin{figure}[ht]
\centering
\includegraphics[scale=0.45]{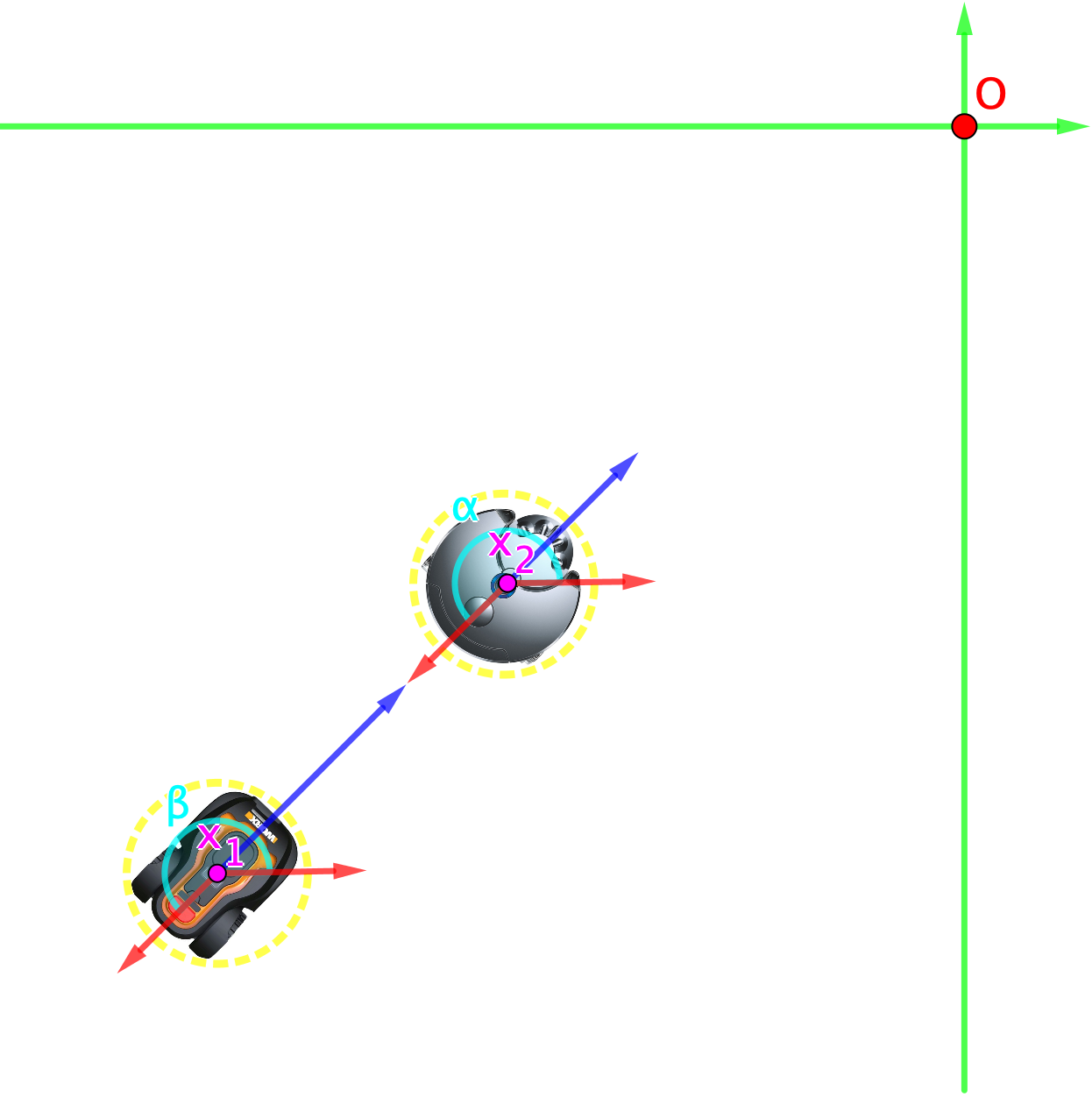}
\caption{Mobile robot model.}
\label{Robot}
\end{figure}

Recall from \cite{cmn1} that the model deals with $n\ge 2$ mobile robots represented as disks of same radius $R$ on the plane. Each robot aims at reaching the target by the shortest pass on the prescribed time interval $[0,T]$ with avoiding (while possibly touching) the other $n-1$ robots, which are treated as obstacles for the robot in question. The model is formalized by considering the configuration vector $x=(x^1,\ldots,x^n)\in\R^{2n}$, where $x^i\in\R^2$ signifies the center of the disk $i$ of each robot with coordinates representation $(\|x^i\|\cos\th_i,\|x^i\|\sin\th_i)$, where $\th_i$ stands for the corresponding constant direction is the smallest positive angle in standard position formed by the positive $x$-axis and vectors $Ox^i$ with $0\in\R^2$ as the target; see Figure~\ref{Robot}. Thus we can describe the trajectory $x^i(t)$ of the $i$-robot by
\begin{eqnarray*}
\ox^i(t)=\big(\|\ox^i(t)\|\cos\th_i,\|\ox^i(t)\|\sin\th_i\big)\;\mbox{ for }\;i=1,\ldots,n,
\end{eqnarray*}
\vspace*{-0.2in}

Following \cite{cmn1}, the dynamic optimization model under consideration can be described in the form of the sweeping control problem $(P)$ with the cost function $\ph(x):=\frac{1}{2}\|x\|^2$ as $x\in\R^{2n}$ and a convex control set $U\subset\R^n$, where---in contrast to \cite{cmn1}---the nonsmooth perturbation mapping $g$ is defined now by
\begin{equation}\label{g}
g(x,u):=\big(s_1|u^1|\cos\th_1,s_1|u^1|\sin\th_1,\ldots,s_n|u^n|\cos\th_n,s_n|u^n|\sin\th_n\big),
\end{equation}
where $s_i$ denotes the speed of robot $i$. Taking into account the convexity of $U$ and the forms of the cost function $\ph$ and of the perturbation mapping $g$ in \eqref{g}, it is easy to conclude by employing the standard variational arguments, which are based on the classical Weierstrass theorem in the weak topology of the space $W^{1,2}([0,1];\R^{2n})\times L^2([0,1];\R^n)$, that problem $(P)$ admits optimal solutions from $W^{1,2}([0,1];\R^{2n})\times L^2([0,1];\R^n)$. Let $(\ox(\cdot),\ou(\cdot))$ be such an optimal solution, which clearly is a relaxed $W^{1,2}\times L^2$-local minimizer of $(P)$ in the sense of Definition~\ref{relaxed}. We are going to apply the results of Theorem~\ref{nsThm6.1*} to determine this solution.\vspace*{-0.03in}

It is not our goal here to investigate the formulated optimal control problem of robotic modeling in generality. Below we confine ourselves to the case study by specifying the initial data of problem $(P)$ as follows:
$n=2,\;x^{01}=\(-30,-30\),\;x^{02}=\(-20,-20\),\;T=6,\;R=12,\;s_1=3,\;s_2=1$, $\th=\alpha=\beta=225^{\circ}$ on Figure~\ref{Robot}, the control set $U$ defined by
\begin{equation*}
U:=\big\{u=(u^1,u^2)\in\R^2\;\big|\;u^1=2u^2,\;-4\le u^1\le 3\big\},
\end{equation*}
and the convex polyhedron $C\subset\R^{4}$ from \eqref{C} given by
\begin{eqnarray*}
C:=\big\{x\in\R^{4}\;\big|\;\la x_*,x\ra\le c\big\}\;\mbox{ with }\;c:=-12\;\mbox{ and }\;x_*:=\(1,1,-1,-1\)\in\R^4.
\end{eqnarray*}
In this case we have $(x^{11}(0)-x^{21}(0))^2+(x^{12}(0)-x^{22}(0))^2=200$ and $t_1>0$ for the first contacting time
\begin{eqnarray*}
t_1:=\min\big\{t\in[0,T]\;\big|\;\|\ox^1(t)-\ox^2(t)\|=2R\big\}.
\end{eqnarray*}
Assume that the robot tends to keep its constant direction and velocity until either touching the other robot (obstacle), or reaching the end of the process at $t=6$.\vspace*{-0.05in}

Applying the necessary optimality conditions of Theorem~\ref{nsThm6.1*}, we arrive at the following relationships:
\begin{itemize}
\item[\bf(1)] $-\(\dot{\ox}^{11}(t),\;\dot{\ox}^{12}(t),\;\dot{\ox}^{21}(t),\;\dot{\ox}^{22}(t)\)=\eta(t)\(1,1,-1,-1\)-\\
\Big(s_1|\ou^1(t)|\cos\th,\;s_1|\ou^1(t)|\sin\th,\;s_2|\ou^2(t)|\cos\th,\;s_2|\ou^2(t)|\sin\th\Big)$ for a.e.\ $t\in[0,6]$.
\item[\bf(2)] $\|\ox^{1}(t)-\ox^{2}(t)\|>2R\Longrightarrow\eta(t)=0$ for a.e.\ $t\in[0,6]$ including $t=6$.
\item[\bf(3)] $\eta(t)>0\Longrightarrow\la x_*,q(t)\ra=c$ for a.e.\ $t\in[0,6)$.
\item[\bf(4)] $\big(\dot{p}(t),-\psi(t)\big)\in\co\partial\big\la q(t),g\big\ra\big(\ox(t),\ou(t)\big),$ for a.e. $t\in[0,6]$, where the set on the right hand-side is calculated as follows
$$\big\{(x(t),\psi^{1}(t),\psi^{2}(t))\in\R^6\big|\;x(t)=0_{\R^4},-s_iq^{i1}(t)\cos\th -s_iq^{i2}(t)\sin\th\le\psi^{i}(t)\le s_iq^{i1}(t)\cos\th +s_iq^{i2}(t)\sin\th\big\},$$
with $i=1,2$, $\th=225^{\circ}$, $x(t)=\(x^{11}(t),x^{12}(t),x^{21}(t),x^{22}(t)\)$, $q(t)=\(q^{11}(t),q^{12}(t),q^{21}(t),q^{22}(t)\)$ and $\psi(t)=\(\psi^{1}(t),\psi^{2}(t)\)$.
\item[\bf(5)] $q(t)=p(t)-\disp\int_{(t,6]}d\gg\(\tau\)x_*$ for all $t\in[0,6]$ except at most a countable subset.
\item[\bf(6)] $\big\la\psi(t),\ou(t)\big\ra=\max_{u\in U}\big\la\psi(t),u\big\ra,\;\mbox{ for a.e. }\;t\in[0,6]$.
\item[\bf(7)] $-p(6)=\lm\ox(6)+\eta(6)x_*$.
\item[\bf(8)] $\eta(6)x_*\in N\big(\ox(6);C)$.
\item[\bf(9)] $\(\lm,p,\|\gg\|\)\ne 0$.
\end{itemize}\vspace*{-0.12in}

Taking into about that the robot directions are constant as well as the assumptions in the model imposed about, we seek for simplicity constant optimal controls. Then follows from the above condition (2) that the function $\eta(\cdot)$ is piecewise constant on $[0,6]$ and admits the representation
\begin{eqnarray*}
\eta(t)=\left\{\begin{array}{ll}
0\quad\mbox{ for }\;t\in[0,t_1),\\
\eta\;\;\mbox{ for }\;t\in[t_1,6].
\end{array}\right.
\end{eqnarray*}
Using now (1), the dynamic equations prior to and after the time $t_1$ can be rewritten as
\begin{eqnarray*}
\left\{\begin{array}{ll}
\dot{\ox}^1(t)=\(s_1|\ou^1|\cos\th,s_1|\ou^1|\sin\th\)\;\mbox{ and}\\
\dot{\ox}^2(t)=\(s_2|\ou^2|\cos\th,s_2|\ou^2|\sin\th\)\;\mbox{ for }\;t\in[0,t_1),
\end{array}\right.
\end{eqnarray*}\vspace*{-0.2in}
\begin{eqnarray*}
\left\{\begin{array}{ll}
\dot{\ox}^1(t)=\(-\eta(t)+s_1|\ou^1|\cos\th,-\eta(t)+s_1|\ou^1|\sin\th\)\;\mbox{ and}\\
\dot{\ox}^2(t)=\(\eta(t)+s_2|\ou^2|\cos\th,\eta(t)+s_2|\ou^2|\sin\th\)\;\mbox{ for }\;t\in[t_1,6].
\end{array}\right.
\end{eqnarray*}
Remembering that the two robots have the same velocities at the first contacting time $t=t_1$ and then maintain them till the end of the process, we get $\dot{\ox}^1(t)=\dot{\ox}^2(t)$ for all $t\in[t_1,6]$. This gives us in turn the following calculation of the corresponding value of $\eta$:
\begin{eqnarray}\label{eq-eta2}
\eta=\left\{\begin{array}{ll}
\frac{1}{2}\(s_1|\ou^1|\cos\th-s_2|\ou^2|\cos\th\)\;&\mbox{ if }\;s_1|\ou^1|\ne  s_2|\ou^2|\;\mbox{ and }\;\cos\th=\sin\th,\\
0&\mbox{ otherwise},
\end{array}\right.
\end{eqnarray}
where the case of $\eta=0$ is trivial. Using the upper formula in \eqref{eq-eta2} and the constant robot and obstacle velocity after the touching moment till reaching the target, we get by the Newton-Leibniz formula that
\begin{eqnarray*}
\left\{\begin{array}{ll}
\ox^1(t)=\(\ox^{11}(0),\ox^{12}(0)\)+\(ts_1|\ou^1|\cos\th,ts_1|\ou^1|\sin\th\)\;\mbox{ and}\\
\ox^2(t)=\(\ox^{21}(0),\ox^{22}(0)\)+\(ts_2|\ou^2|\cos\th,ts_2|\ou^2|\sin\th\)\;\mbox{ for }\;t\in[0,t_1),
\end{array}\right.
\end{eqnarray*}\vspace*{-0.2in}
\begin{eqnarray*}
\left\{\begin{array}{ll}
\ox^1(t)=\(\ox^{11}(0),\ox^{12}(0)\)+\big(ts_1|\ou^1|\cos\th-\eta(t-t_1),ts_1|\ou^1|\sin\th-\eta(t-t_1)\big)\;\mbox{ and}\\
\ox^2(t)=\(\ox^{21}(0),\ox^{22}(0)\)+\big(ts_2|\ou^2|\cos\th+\eta(t-t_1),ts_2|\ou^2|\sin\th+\eta(t-t_1)\big)\;\mbox{ for }\;t\in[t_1,6].
\end{array}\right.
\end{eqnarray*}
Recalling that $\|\ox^2(t_1)-\ox^1(t_1)\|=2R$ and taking into account that we are considering the sum norm, the latter leads us to the following quadratic equation for $t_1$:
\begin{eqnarray}\label{t_1}
\begin{array}{ll}
|\ox^{21}(0)-\ox^{11}(0)+t_1\(s_2|\ou^2|-s_1|\ou^1|\)\cos\th |+
|\ox^{22}(0)-\ox^{12}(0)+t_1\(s_2|\ou^2|-s_1|\ou^1|\)\sin\th |=2R,
\end{array}
\end{eqnarray}
which makes the connection between $t_1$ and the control $\ou=(\ou^1,\ou^2)$ via the given model data. It follows from (4) and (6) that the control $\(\ou_1,\ou_2\)$ under consideration can be either $(3,3/2)$ or $(-4,-2)$.

{\bf Case 1:} We deduce from \eqref{t_1} that $t_1\approx 0.38$. Moreover, we get from \eqref{eq-eta2} that
\begin{eqnarray*}
\eta=\frac{1}{2}\Big(3|\ou^1|\Big(-\frac{\sqrt{2}}{2}\Big)-|\ou^2|\Big(-\frac{\sqrt{2}}{2}\Big)\Big)=-\frac{5\sqrt{2}}{4}|\ou^2|\ne 0,
\end{eqnarray*}
which gives us $\eta=-\dfrac{15\sqrt{2}}{8}$. The above calculations bring us to the expressions:
\begin{eqnarray*}
\left\{\begin{array}{ll}
(\ou^1,\ou^2)&=\(3,1.5\),\\
\ox^1(t)&\approx\(-30-6.36t,-30-6.36t\),\quad t\in [0,0.38),\\
\ox^1(t)&\approx\(-31.01-3.71t,-31.01-3.71t\),\quad t\in [0.38,6],\\
\ox^2(t)&\approx\(-20-1.06t,-20-1.06t\),\quad t\in [0,0.38),\\
\ox^2(t)&\approx\(-18.99-3.71t,-18.99-3.71t\),\quad t\in[0.38,6].
\end{array}\right.
\end{eqnarray*}\vspace*{-0.1in}

{\bf Case 2:} It follows from \eqref{t_1} that $t_1\approx 0.28$. Moreover, we obtain from \eqref{eq-eta2} that $\eta=-\dfrac{5\sqrt{2}}{2}$. The above calculations bring us to the expressions:
\begin{eqnarray*}
\left\{\begin{array}{ll}
(\ou^1,\ou^2)&=\(-4,-2\),\\
\ox^1(t)&\approx\(-30-8.49t,-30-8.49t\),\quad t\in [0,0.28),\\
\ox^1(t)&\approx\(-30.99-4.95t,-30.99-4.95t\),\quad t\in [0.28,6],\\
\ox^2(t)&\approx\(-20-1.41t,-20-1.41t\),\quad t\in [0,0.28),\\
\ox^2(t)&\approx\(-19.01-4.95t,-19.01-4.95t\),\quad t\in[0.28,6].
\end{array}\right.
\end{eqnarray*}\vspace*{-0.1in}

Comparing the cost function in the above two cases, we obtain that $(\ou^1,\ou^2)=(3,1.5)$ is the optimal control to this problem. To make a conclusion about the optimality of the obtained solution, we have to check the fulfillment of all the other necessary conditions of Theorem~\ref{nsThm6.1*}. It follows that the corresponding adjoint arc $q(\cdot)$ can be calculated from (4) and (6) with choosing $\psi(t)=\ou(t)$. This gives us the values $q^{11}(t)=-1,\;q^{12}(t)=-1,\;q^{21}(t)=-\sqrt{2}$, and $q^{22}(t)=-\sqrt{2}$.

Then we deduce from (4) and (7) that $\dot{p}(t)=0\;\mbox{ for a.e. }t\in[0,T]$ and $p(6)=-\lm\ox(6)-\eta x_*$ with $\eta=-\dfrac{15\sqrt{2}}{8}$ and $x_*=\(1,1,-1,-1\)$. Hence equation (5) reads as
\begin{eqnarray*}
\gamma([t,6])=p(t)-q(t)\;\mbox{ a.e. on }\;[0,6].
\end{eqnarray*}
Combining the latter with the above calculations tells us that
\begin{eqnarray*}
\gamma([t,6])&\approx&\(55.92,\;55.92,\;38.60,\;38.60\)-\(-1,\;-1,\;-1.41,\;-1.41\)\\
&\approx&\(56.92,\;56.92,\;40.01,\;40.01\)
\end{eqnarray*}
for $0.38\le t\le 6$ and $\lambda=1$. It confirms that the calculated motion hits the boundary of the state constraint at the time $t_1\approx 0.38$ and stays there until the end of the process. Summarizing all the above, we conclude that Theorem~\ref{nsThm6.1*} allows us, under the assumptions made, to single out a feasible pair $(\ox(\cdot),\ou(\cdot))$ satisfying the obtained necessary optimality conditions. Taking into account the existence of optimal solutions to this problems ensures that the obtained pair $(\ox(\cdot),\ou(\cdot))$ is the one.

{\bf Acknowledgements.} The authors are gratefully indebted to Tan Cao and Giovanni Colombo for many useful discussions on the material presented in this paper.
\end{document}